\documentclass{amsart}
\usepackage{amsfonts}
\usepackage{amssymb}
 \newtheorem{cor}{Corollary}
 \newtheorem{theorem}{Theorem}
 \newtheorem{lemma}{Lemma}

 \newtheorem{remark}{Remark}
 
 \newcommand{\bee}[1]{\begin{equation}\label{#1}}
 \newcommand{\ene}{\end{equation}}

\begin{document}

\title{*-Graded Capelli Polynomials  and their Asymptotics }
\author[Benanti]{F. S. Benanti}\address{Dipartimento di Matematica e Informatica
\\ Universit\`a di Palermo \\
via Archirafi, 34\\ 90123
Palermo,Italy}
\email{francescasaviella.benanti@unipa.it}

\author[Valenti]{A. Valenti }\address{Dipartimento di Ingegneria  \\ Universit\`a di Palermo
\\Viale delle Scienze \\
90128 Palermo, Italy}
 \email{angela.valenti@unipa.it}


\date{}

\begin{abstract}
Let  $F\langle Y \cup Z, \ast
\rangle$ be the free $\ast$-superalgebra over a
field $F$ of characteristic zero and let $ \Gamma^\ast_{M^{\pm},
L^{\pm}} $ be the $T^\ast_{\mathbb{Z}_2}$-ideal generated by the set of
the $\ast$-graded Capelli polynomials
 $Cap^{(\mathbb{Z}_2, \ast)}_{M^+}
[Y^+,X]$, $Cap^{(\mathbb{Z}_2, \ast)}_{M^-} [Y^-,X]$,
$Cap^{(\mathbb{Z}_2, \ast)}_{L^+} [Z^+,X]$, $Cap^{(\mathbb{Z}_2,
\ast)}_{L^-} [Z^-,X]$ alternating on $M^+$ symmetric variables
of homogeneous degree zero, on $M^-$ skew variables of homogeneous
degree zero, on $L^+$ symmetric variables of homogeneous degree
one and on $L^-$ skew variables of homogeneous degree one,
respectively.
 We study the asymptotic behavior of the
sequence of $\ast$-graded codimensions of $\Gamma^\ast_{M^{\pm},
L^{\pm}}.$ In particular we prove that the $\ast$-graded codimensions of
the finite dimensional simple $\ast$-superalgebras are asymptotically
equal to the $\ast$-graded codimensions of $\Gamma^\ast_{M^{\pm},
L^{\pm}}$, for some fixed natural
numbers $M^+, M^-, L^+$ and $L^-$.

\end{abstract}

\maketitle

\section{Introduction}
This paper is devoted to the study of the $\ast$-superalgebras, i.e.
superalgebras endowed with a graded involution, and the asymptotic
behavior of their $\ast$-graded codimensions. If $A$ is an algebra
over a field $F$ of characteristic zero an effective way of
measuring the polynomial identities satisfied by $A$ is provided by
its sequence of codimensions $\{c_n(A)\}_{n\geq 1}$ whose $n$-th
therm is the dimension of the space of multilinear polynomials in
$n$ variables in the corresponding relatively free algebra of
countable rank. Such sequence was introduced by Regev in \cite{Re1}
and, in characteristic zero, gives a quantitative measure of the
identities satisfied by a given algebra. The most important result
of the sequence of codimensions proved in \cite{Re1} states that if
$A$ is a PI-algebra, i.e. it satisfies a non trivial polynomial
identity, then $\{c_n(A)\}_{n\geq 1}$ is exponential bounded.
  Later,  Giambruno
and Zaicev (\cite{GZ00}, \cite{GZ0})  answered in a positive way to a well known  conjecture of Amitsur proving the existence and the integrality of
$$
\mathrm{exp}(A)=\displaystyle{\lim_{n\rightarrow \infty}\root n
\of {c_n(A)}}
$$
the {exponent} of $A$. These results, in the last years, have been
extended to algebras with an additional structure as algebras with
involution (\cite{AG}, \cite{GPV}), superalgebras (\cite{BGP}) and
more generally algebras graded by a group ( \cite{AGL}, \cite{DN},
\cite{GL}, \cite{GZ1} ), algebras with a generalised $H$-action
(\cite{Go}), superalgebras with graded involution (\cite{S}) and
superalgebras with superinvolution (\cite{I}).

Let $A= A_{0}^+\oplus  A_{0}^- \oplus A_{1}^+\oplus  A_{1}^-$ be a
$\ast$-superalgebra and let $c_n^{(\mathbb{Z}_2,\ast)}(A)$,
$n=1,2,\ldots$, be its sequence of $\ast$-graded codimensions. If
$A$ is a PI-algebra it can be easily proved that the relation
between codimensions and $\ast$-graded codimensions is given by
$c_n(A)\leq c_n^{(\mathbb{Z}_2,\ast)}(A) \leq 4^n c_n(A)$. Hence, as
in the ordinary case, the sequence of $\ast$-graded codimensions is
exponentially bounded. Moreover, since a $\ast$-superalgebra can be
viewed as an algebra with a generalized FG-action where
$G=\mathbb{Z}_2 \times \mathbb{Z}_2$ acts on it by automorphism and
antiautomorphism, in the finite dimensional case, the existence of
the $\ast$-graded exponent has been confirmed by Gordienko in
\cite{Go}.

Let $M^+$, $M^-$, $L^+$ and $L^-$ be natural numbers  and let's
denote by $\Gamma^\ast_{M^{\pm},
L^{\pm}}$ the
$T^\ast_{\mathbb{Z}_2}$-ideal generated by the set of the $\ast$-graded
Capelli polynomials
 $Cap^{(\mathbb{Z}_2, \ast)}_{M^+}
[Y^+,X]$, $Cap^{(\mathbb{Z}_2, \ast)}_{M^-} [Y^-,X]$,
$Cap^{(\mathbb{Z}_2, \ast)}_{L^+} [Z^+,X]$, $Cap^{(\mathbb{Z}_2,
\ast)}_{L^-} [Z^-,X]$ alternating on $M^+$ symmetric variables
of homogeneous degree zero, on $M^-$ skew variables of homogeneous
degree zero, on $L^+$ symmetric variables of homogeneous degree
one and on $L^-$ skew variables of homogeneous degree one, respectively. In this
paper we find a relation among the $\ast$-graded codimensions of the
finite dimensional simple $\ast$-superalgebras and  the
$\ast$-graded codimensions of $\Gamma^\ast_{M^{\pm},
L^{\pm}}$ proving their asymptotic equality.
 Recall that two sequences $a_n$, $b_n$, $n=1,2,\ldots$, are
asymptotically equal, $a_n \simeq b_n$, if $ \lim_{n\rightarrow
+\infty}\frac{a_n}{b_n}=1. $ In the ordinary case
(see \cite{GZ2}) it was proved the asymptotic equality between the
codimensions of the Capelli polynomials $Cap_{k^2+1}$ and the
codimensions of the matrix algebra $M_k(F).$ In \cite{Be} this
result was extended to finite dimensional simple superalgebras
and in \cite{BV} the authors found  similar result in the case of algebras with involution (for a survey see \cite{BV2}). The link between the
asymptotic of the codimensions of the Amitsur's Capelli-type
polynomials and the verbally prime algebras was studied in
\cite{BS}.

\section{Preliminaries}

Throughout this paper, $F$ will be a field of characteristic zero
and $A$ an associative algebra over $F$. We say that $A$ is a
\textit{$\mathbb{Z}_2$-graded algebra} or a \textit{superalgebra} if
it can be decomposed into a direct sum of subspaces $A=A_{0}\oplus
A_{1}$ such that $A_{0} A_{0}+ A_{1} A_{1}\subseteq A_{0}$ and
$A_{0} A_{1}+ A_{1} A_{0}\subseteq A_{1}$. The elements of $A_{0}$
are called \textit{ homogeneous of degree zero} (\textit{even
elements}) and those of $A_{1}$ \textit{homogeneous of degree one}
(\textit{odd elements}).

Recall that an \textit{involution} $\ast$ on an algebra $A$ is just an antiautomorphism on $A$ of order at most $2$. We write $A^+=\{a \in A \,\, | \,\, a^\ast =a\}$ and
$A^-=\{a \in A \,\, | \,\, a^\ast =-a\}$ for the set of \textit{symmetric} and \textit{skew symmetric} elements of $A$ respectively.

Given a superalgebra $A=A_{0}\oplus A_{1}$ endowed with an
involution $\ast$, we say that $\ast$ is a \textit{graded
involution} if it preserves the homogeneous components of $A$, i.e.
if $A_{i}^\ast \subseteq A_{i}$, $i=0,1$. A superalgebra endowed
with a graded involution is called \textit{$\ast$-superalgebra}. It
is clear that a superalgebra $A$ is a $\ast$-superalgebra if and
only if the subspaces $A^+$ and $A^-$ are graded subspaces, i.e.
$A^+=A^+_{0}\oplus A^+_{1}$ and $A^-=A^-_{0}\oplus A^-_{1}$. Thus,
since char $F$= 0, the $\ast$-superalgebra $A$ can be written as
$$
A=A_{0}^+\oplus A_{0}^-\oplus A_{1}^+\oplus A_{1}^-
$$
where, for  $i=0,1$, $A_{i}^+=\{a \in A_{i} \,\, | \,\, a^\ast =a\}$ and $A_{i}^-=\{a \in A_{i} \,\, | \,\, a^\ast =-a\}$ denote the sets of homogeneous symmetric and
skew elements of $A_{i}$, respectively.
We remark that an algebra with involution $\ast$ and trivial $\mathbb{Z}_2$-grading is a $\ast$-superalgebra.

Let $A$ be a $\ast$-superalgebra and let $I$ be an ideal of $A$, we  say that $I$ is a \textit{$\ast$-graded ideal }of $A$ if it is homogeneous in the $\mathbb{Z}_2$-grading
and invariant under $\ast$. Moreover  $A$ is called \textit{simple} $\ast$-superalgebra if $A^2\neq \{0\}$ and it has no non-zero $\ast$-graded ideals.

Let $X=\{x_1,x_2,\ldots\}$ be
a countable set of non commutative variables and $F\langle X\rangle$ the free associative algebra on $X$ over $F$. We write $X=Y \cup Z$ as the disjoint union of two
countable
sets of variables $Y=\{y_1,y_2,\ldots\}$ and $Z=\{z_1,z_2,\ldots\}$, then
$F\langle X\rangle=F\langle Y \cup Z \rangle=\langle y_1,z_1,y_2,z_2,\ldots\rangle$ has a natural
structure of free superalgebra if we require that the variables
from $Y$ have degree zero and the variables from $Z$ have degree
one. This algebra is said to be the \textit{free superalgebra} over $F$.
Moreover, if we write each set as the disjoint union of two other infinite sets of symmetric and skew elements, respectively, then we obtain the \textit{free
$\ast$-superalgebra} $$
F\langle Y \cup Z, \ast \rangle= F\langle y_1^+,y_1^-,z_1^+, z_1^-, \ldots  \rangle
$$
where $y_i^+=y_i+y_i^\ast$ denotes a symmetric variable of even degree, $y_i^-=y_i-y_i^\ast$ a skew variable of even degree, $z_i^+=z_i+z_i^\ast$ a symmetric variable of odd
degree and  $z_i^-=z_i-z_i^\ast$ a skew variable of
odd degree.

 An element  $f=f(y_1^+, \ldots,y_n^+,y_1^-,
\ldots,y_m^-,z_1^+, \ldots, z_p^+,z_1^-, \ldots, z_q^-)$ of
$F\langle Y \cup Z, \ast \rangle$ is a \textit{$\ast$-graded
polynomial identity }for a $\ast$-superalgebra $A$ if
$$
f(a_{1,0}^+,
\ldots , a_{n,0}^+,a_{1,0}^-, \ldots, a_{m,0}^-,a_{1,1}^+, \ldots ,
a_{p,1}^+, a_{1,1}^-, \ldots, a_{q,1}^-) = 0_A
$$

\noindent  for every $a_{1,0}^+, \ldots , a_{n,0}^+\in A_{0}^+$,
$a_{1,0}^-, \ldots, a_{m,0}^- \in A_{0}^-$, $a_{1,1}^+, \ldots ,
a_{p,1}^+\in A_{1}^+$,  $a_{1,1}^-, \ldots, a_{q,1}^- \in A_{1}^-$
and we write $f\equiv 0$. The set of all  $\ast$-graded polynomial
identities satisfied by $A$
$$
Id^\ast_{\mathbb{Z}_2}(A)=\{ f \in F\langle Y \cup Z, \ast \rangle \,\, | \,\, f\equiv 0 \,\, \mathrm{on} \, A \}
$$

\noindent is an ideal of $F\langle Y \cup Z, \ast \rangle$ called
the \textit{ideal of $\ast$-graded  identities of $A$}. It is easy
to show that $Id^\ast_{\mathbb{Z}_2}(A)$ is a
\textit{$T^\ast_{\mathbb{Z}_2}$-ideal} of $F\langle Y \cup Z, \ast
\rangle$, i.e. a two-sided ideal invariant under all endomorphisms
of the free $\ast$-superalgebra that preserve the superstructure and
commute with the graded involution $\ast$. Now, let
$$
P^{(\mathbb{Z}_2, \ast)}_n=\{w_{\sigma(1)}, \ldots , w_{\sigma(n)} \,\, | \,\, \sigma \in S_n, \, w_i \in \{y_i^+,y_i^-,z_i^+,z_i^-\}, \, i=1, \ldots , n\}
$$

\noindent be the space of multilinear polynomials of degree $n$  in
the variables
 $y_1^+$, $y_1^-$, $z_1^+$, $z_1^-$,$ \ldots,$\, $y_n^+$, $y_n^-$, $z_n^+$, $z_n^- $,  (i.e., $y_i^+$, $y_i^-$, $z_i^+$  or $z_i^-$ appears in each
monomial at degree $1$). Since char $F$= 0, it is well known that
$Id^\ast_{\mathbb{Z}_2}(A)$ is completely determined by its
multilinear polynomials, then the study of
$Id^\ast_{\mathbb{Z}_2}(A)$ is equivalent to that of
$Id^\ast_{\mathbb{Z}_2}(A)\cap P^{(\mathbb{Z}_2, \ast)}_n$ for all
$n\geq 1$. As in the ordinary case (see \cite{Re1}), one defines the
\textit{$n$-th $\ast$-graded codimension} $c_n^{(\mathbb{Z}_2,
\ast)}(A)$ of the $\ast$-superalgebra $A$ as
$$
c_n^{(\mathbb{Z}_2, \ast)}(A)= \dim_F  \frac{P^{(\mathbb{Z}_2,
\ast)}_n }{P^{(\mathbb{Z}_2, \ast)}_n \cap
Id^\ast_{\mathbb{Z}_2}(A)}.
$$

If $A$ is a PI-algebra, i.e. satisfies an ordinary polynomial
identity, then the sequence $\{c_n^{(\mathbb{Z}_2,
\ast)}(A)\}_{n\geq 1}$ is exponentially bounded (see \cite[Lemma
3.1]{GSV}). If $A$ is a finite dimensional PI-algebra,  Gordienko in
\cite{Go} proved that

$$
\mathrm{exp}_{\mathbb{Z}_2}^\ast(A)=\displaystyle{\lim_{n\rightarrow \infty}\root n
\of{c_n^{(\mathbb{Z}_2, \ast)}(A)}}
$$

\medskip
\noindent exists and is a non-negative integer which is called the
\textit{$\ast$-graded exponent}
 of the $\ast$-superalgebra $A$.
 It is often more useful to study $\ast$-superalgebras up to $\ast$-graded PI-equivalence, then it is convenient to use the language of varieties. Let $I$ be a
 $T_{\mathbb{Z}_2}^\ast$-ideal
 of $F\langle Y \cup Z , \ast \rangle$ and $\mathcal{V}_{\mathbb{Z}_2}^\ast $ the {\textit{variety of $\ast$-superalgebras}
 associated to $I$, i.e. the class of all the $\ast$-superalgebras $A$ such that $I$ is contained in $Id_{\mathbb{Z}_2}^\ast(A)$. We put
 $I=Id_{\mathbb{Z}_2}^\ast(\mathcal{V}_{\mathbb{Z}_2}^\ast)$.
 When $Id_{\mathbb{Z}_2}^\ast(\mathcal{V}_{\mathbb{Z}_2}^\ast)=Id_{\mathbb{Z}_2}^\ast(A)$ we say that the variety  $\mathcal{V}_{\mathbb{Z}_2}^\ast$  is generated by the
 $\ast$-superalgebra $A$ and we write $\mathcal{V}_{\mathbb{Z}_2}^\ast=\mathrm{var}_{\mathbb{Z}_2}^\ast (A)$ and set
 $\mathrm{exp}_{\mathbb{Z}_2}^\ast(\mathcal{V}_{\mathbb{Z}_2}^\ast)=\mathrm{exp}_{\mathbb{Z}_2}^\ast(A)$ the $\ast$-graded exponent of the variety
 $\mathcal{V}_{\mathbb{Z}_2}^\ast$, if $\mathrm{exp}_{\mathbb{Z}_2}^\ast(A)$ exists.

Now, if $f\in F\langle Y \cup
Z , \ast\rangle$ we denote by $\langle f\rangle_{\mathbb{Z}_2}^\ast$ the $T_{\mathbb{Z}_2}^\ast$-ideal
generated by $f$. Also for a set of polynomials $V\subset F\langle
Y \cup Z, \ast \rangle$ we write $\langle V\rangle_{\mathbb{Z}_2}^\ast$ to indicate
the $T_{\mathbb{Z}_2}^\ast$-ideal generated by $V$.

In PI-theory a prominent role is played by the Capelli polynomial. Let us recall that, for any positive integer $m$, the \textit{$m$-th Capelli polynomial} is the element of
$F\langle X\rangle$ defined as
$$
Cap_m [T,X]= Cap_m (t_1, \ldots , t_m ; x_1, \ldots , x_{m-1}) =$$
$$ =
\sum_{\sigma \in S_m} \mathrm{(sgn
\sigma)}t_{\sigma(1)}x_1t_{\sigma(2)}\cdots t_{\sigma(m-1)}
x_{m-1}t_{\sigma(m)}
$$
where $S_m$ is the symmetric group on $\{1, \ldots , m\}$.
In particular we write

$$
Cap_m^{(\mathbb{Z}_2, \ast)} [Y^+,X],\,\,\,  Cap_m^{(\mathbb{Z}_2, \ast)}[Y^-,X], \,\,\, Cap_m^{(\mathbb{Z}_2, \ast)} [Z^+,X]\,\,\,\textrm{ and} \,\,\, Cap_m^{(\mathbb{Z}_2,
\ast)}[Z^-,X]
$$

\noindent to indicate the \textit{$m$-th $\ast$-graded Capelli
polynomial} alternating in the
 symmetric variables of degree zero $y_1^+, \ldots , y_m^+$, in the skew variables of degree zero $y_1^-, \ldots , y_m^-$, in the symmetric variables of degree one $z_1^+,
 \ldots , z_m^+$
 and in the skew variables of degree one  $z_1^-, \ldots , z_m^-$, respectively ($x_1, \ldots ,x_{m-1}$ are arbitrary variables).
 Let $\overline{Cap}_m^{(\mathbb{Z}_2, \ast)} [Y^+,X]$
 denote the set of $2^{m-1}$ polynomials obtained from $Cap_m^{(\mathbb{Z}_2, \ast)} [Y^+,X]$ by deleting any subset of variables $x_i$  (by
evaluating the variables $x_i$ to $1$ in all possible way).  In a
similar way we define $\overline{Cap}_m^{(\mathbb{Z}_2, \ast)}[Y^-,X]$,
$\overline{Cap}_m^{(\mathbb{Z}_2, \ast)} [Z^+,X]$ and
$\overline{Cap}_m^{(\mathbb{Z}_2, \ast)}[Z^-,X]$. If $M^+$, $M^-$, $L^+$
and $L^-$ are natural numbers,  we denote by

$$
\Gamma^\ast_{M^{\pm},
L^{\pm}}=\langle \overline{Cap}_{M^+}^{(\mathbb{Z}_2, \ast)} [Y^+,X],  \overline{Cap}_{M^-}^{(\mathbb{Z}_2, \ast)}[Y^-,X],   \overline{Cap}_{L^+}^{(\mathbb{Z}_2,
\ast)} [Z^+,X], \overline{Cap}_{L^-}^{(\mathbb{Z}_2, \ast)}[Z^-,X]\rangle_{\mathbb{Z}_2}^\ast
$$

\medskip
\noindent the $T_{\mathbb{Z}_2}^\ast$-ideal
generated by  $\overline{Cap}_{M^+}^{(\mathbb{Z}_2, \ast)} [Y^+,X]$,  $\overline{Cap}_{M^-}^{(\mathbb{Z}_2, \ast)}[Y^-,X]$,   $\overline{Cap}_{L^+}^{(\mathbb{Z}_2, \ast)}
[Z^+,X]$  and \\ $\overline{Cap}_{L^-}^{(\mathbb{Z}_2, \ast)}[Z^-,X]$.

The purpose of this paper is to find a close relation among the
asymptotic behavior of the $\ast$-graded codimensions of any finite
dimensional simple $\ast$-superalgebra $A=A_{0}^+\oplus
A_{0}^-\oplus A_{1}^+\oplus A_{1}^-$ and the asymptotic behavior of
the $\ast$-graded codimensions of $\Gamma^\ast_{M^{\pm}+1,
L^{\pm}+1}$, where $M^+=\mathrm{dim}_FA_{0}^+$,
$M^-=\mathrm{dim}_FA_{0}^-$, $L^+=\mathrm{dim}_FA_{1}^+$ and
$L^-=\mathrm{dim}_FA_{1}^-$. More precisely,  we characterize the
$T_{\mathbb{Z}_2}^\ast$-ideal $Id_{\mathbb{Z}_2}^\ast(A)$ showing
that

$$
\Gamma^\ast_{M^{\pm}+1,
L^{\pm}+1}=Id_{\mathbb{Z}_2}^\ast(A\oplus D),
$$
where $D$ is a finite dimensional $\ast$-superalgebra such that $\mathrm{exp}_{\mathbb{Z}_2}^\ast(D)<\mathrm{exp}_{\mathbb{Z}_2}^\ast(A)$. Moreover we obtain the asymptotic
equality

$$
c_n^{(\mathbb{Z}_2, \ast)}(\Gamma^\ast_{M^{\pm}+1,
L^{\pm}+1})\simeq c_n^{(\mathbb{Z}_2, \ast)}(A).
$$

 \bigskip

\section{Basic Results}

Let $A$ be a finite dimensional  $\ast$-superalgebra  over a field
$F$ of characteristic zero. From now on we assume that $F$ is
algebraically closed. In fact, since $F$ has characteristic zero,
$Id_{\mathbb{Z}_2}^\ast(A)=Id_{\mathbb{Z}_2}^\ast(A\otimes_F L)$ for
any extension field $L$ of $F$ then also the $\ast$-graded
codimensions of $A$ do not change upon extension of the base field.
By the generalization of the Wedderburn-Malcev Theorem (see
\cite[Theorem 7.3]{GSV}), we can write $A=A_1 \oplus \cdots \oplus
A_s+J$, where $A_1, \ldots, A_s$ are simple $\ast$-superalgebras and
$J=J(A)$ is the Jacobson radical of $A$ which is a $\ast$-graded
ideal.

 We say that a subalgebra  $A_{i_1}\oplus \cdots \oplus  A_{i_k}$ of $A$, where $A_{i_1}, \ldots ,   A_{i_k}$ are distinct simple
components, is \textit{admissible} if for some permutation $(l_1,
\ldots , l_k)$ of $(i_1, \ldots , i_k)$ we have that $A_{l_1}J
\cdots J A_{l_k}\neq 0$. Moreover, if $A_{i_1}\oplus \cdots \oplus
A_{i_k}$ is an admissible subalgebra of $A$ then $A'=A_{i_1}\oplus
\cdots \oplus  A_{i_k}+J$ is called a \textit{reduced} algebra.

 The
notion of admissible $\ast$-superalgebra is closely linked to that
of $\ast$-graded exponent in fact, in \cite{Go}, Gordienko proved
that $\mathrm{exp}_{\mathbb{Z}_2}^\ast(A)=d$ where $d$ is the
maximal dimension of an admissible subalgebra of $A$.
It follows immediately that

\smallskip

 \begin{remark} If $A$ is a simple $*$-superalgebra then
$\mathrm{exp}_{\mathbb{Z}_2}^\ast(A) = \mathrm{dim}_F A$.
\end{remark}

\smallskip

By \cite[Theorem 5.3]{GIL2} the Gordienko's result on the
existence of the $\ast$-graded exponent can be actually extended to
any finitely generated PI-$\ast$-superalgebra since it satisfies the
same $\ast$-graded polynomial identities of a finite-dimensional
$\ast$-superalgebra.

In  \cite{GZ2} it was showed that reduced superalgebras are building
blocks of any proper variety. Here we obtain the analogous result
for varieties of $\ast$-superalgebras.

 Let's first start with the following

\smallskip

\begin{lemma} \label{codimensioni}
Let $A$ and $B$ be  $\ast$-superalgebras satisfying an ordinary
polynomial identity. Then
$$
c_n^{(\mathbb{Z}_2,\ast)}(A), c_n^{(\mathbb{Z}_2,\ast)}(B) \leq
c_n^{(\mathbb{Z}_2,\ast)}(A\oplus B) \leq
c_n^{(\mathbb{Z}_2,\ast)}(A)+c_n^{(\mathbb{Z}_2,\ast)}(B).
$$
If $A$ and $B$ are finitely generated $\ast$-superalgebras, then
$$
\mathrm{exp}_{\mathbb{Z}_2}^\ast(A\oplus
B)=max\{\mathrm{exp}_{\mathbb{Z}_2}^\ast(A),
\mathrm{exp}_{\mathbb{Z}_2}^\ast(B)\}.
$$

\end{lemma}
\noindent {\bf Proof.} The proof  is the same of the proof of the  Lemma 1 in \cite{GZ2}.

\bigskip

We have the following

\smallskip

\begin{theorem} \label{decomposizione}
Let $A$ be a finitely generated  $\ast$-superalgebra satisfying an
ordinary polynomial identity. Then there exists a finite number of
reduced $\ast$-superalgebras $B_1, \ldots, B_t$ and a finite
dimensional $\ast$-superalgebra $D$ such that
$$
\mathrm{var}_{\mathbb{Z}_2}^\ast(A)=\mathrm{var}_{\mathbb{Z}_2}^\ast (B_1\oplus \cdots \oplus B_t \oplus D)
$$
\noindent with
$\mathrm{exp}_{\mathbb{Z}_2}^\ast(A)=\mathrm{exp}_{\mathbb{Z}_2}^\ast(B_1)=\cdots=\mathrm{exp}_{\mathbb{Z}_2}^\ast(B_t)$
and
$\mathrm{exp}_{\mathbb{Z}_2}^\ast(D)<\mathrm{exp}_{\mathbb{Z}_2}^\ast(A)$.
\end{theorem}

\noindent {\bf Proof.} The proof follows closely the proof given in
\cite[Theorem 3]{Be}.  Since $A$ is a finitely generated
$\ast$-superalgebra,  by \cite{GIL2}, there exists a finite
dimensional $\ast$-superalgebra $B$ such that
$Id_{Z\mathbb{}_2}^{\ast}(A) =  Id_{Z\mathbb{}_2}^{\ast}(B)$.
Therefore we may assume that $A= A_{0}^+\oplus  A_{0}^- \oplus
A_{1}^+\oplus  A_{1}^-$ is a finite dimensional $\ast$-superalgebra
over $F$ satisfying an ordinary polynomial identity. Also, by
\cite[Theorem 7.3]{GSV} we can write
$$
A=A_1 \oplus \cdots \oplus A_s+J
$$
where $A_1, \ldots A_s$ are simple $\ast$-superalgebras and $J=J(A)$ is the Jacobson radical of $A$ which is a $\ast$-graded ideal.
Let $\mathrm{exp}_{\mathbb{Z}_2}^\ast(A)=d$.  Then there
exist distinct
simple $\ast$-superalgebras $A_{j_1}, \ldots A_{j_k} $ such that
$$
A_{j_1}J \cdots JA_{j_k}\neq 0 \,\,\,\,\,\,\,\mathrm{ and} \,\,\,\,\,\,\, \mathrm{dim}_F(A_{j_1}\oplus \cdots \oplus A_{j_k})=d.
$$
Let $\Gamma_1, \ldots , \Gamma_t$ be all possible subset of $\{1,
\ldots , s\}$ such that, if  $\Gamma_j=\{j_1, \ldots ,j_k\}$ then
$\mathrm{dim}_F(A_{j_1}\oplus \cdots \oplus A_{j_k})=d$ and
$A_{\sigma(j_1)}J \cdots JA_{\sigma(j_k)}\neq 0$ for some
permutation $\sigma \in S_k$. For any such $\Gamma_j$, $j=1, \ldots
, t,$ then we put $B_j=A_{j_1}\oplus \cdots \oplus A_{j_k} + J$. It
follows, by the characterization of the $\ast$-graded exponent, that
$$
\mathrm{exp}_{\mathbb{Z}_2}^\ast(B_1)=\cdots =\mathrm{exp}_{\mathbb{Z}_2}^\ast(B_t)=d=\mathrm{exp}_{\mathbb{Z}_2}^\ast(A).
$$
Let $D=D_1 \oplus \cdots \oplus D_p$, where $D_1, \ldots , D_p$ are
all $\ast$-graded subalgebras  of $A$ of the type $A_{i_1}\oplus
\cdots \oplus A_{i_r} + J$, with $1\leq i_1 <\cdots < i_r \leq s$
and $\mathrm{dim}_F(A_{i_1}\oplus \cdots \oplus A_{i_r})<d$. Then,
by the previous lemma, we have that
$\mathrm{exp}_{\mathbb{Z}_2}^\ast(D)<\mathrm{exp}_{\mathbb{Z}_2}^\ast(A)$.
Now, we want to  prove that $\mathrm{exp}_{\mathbb{Z}_2}^\ast(B_1
\oplus \cdots \oplus B_t\oplus
D)=\mathrm{exp}_{\mathbb{Z}_2}^\ast(A)$. The inclusion
$$
\mathrm{var}_{\mathbb{Z}_2}^\ast(B_1 \oplus \cdots \oplus B_t\oplus D)\subseteq \mathrm{var}_{\mathbb{Z}_2}^\ast(A)
$$
follows since $D, B_i\in \mathrm{var}_{\mathbb{Z}_2}^\ast(A)$, $\forall i=1, \ldots , t$.

\noindent Let's consider a multilinear polynomial $f=f(y_1^+,
\ldots,y_n^+,y_1^-, \ldots,y_m^-,z_1^+, \ldots, z_p^+,z_1^-, \ldots,
z_q^-)$   such that $f\not\in Id_{\mathbb{Z}_2}^\ast(A)$. We shall
prove that $f \not\in Id_{\mathbb{Z}_2}^\ast(B_1 \oplus \cdots
\oplus B_t\oplus D) $. Since $f\not\in Id_{\mathbb{Z}_2}^\ast(A)$
there exist $a_{1,0}^+, \ldots , a_{n,0}^+\in A_{0}^+$,
$a_{1,0}^-, \ldots, a_{m,0}^- \in A_{0}^-$,  $a_{1,1}^+, \ldots ,
a_{p,1}^+\in A_{1}^+$ and $a_{1,1}^-, \ldots, a_{q,1}^- \in
A_{1}^-$ such that
$$
f(a_{1,0}^+, \ldots , a_{n,0}^+,a_{1,0}^-, \ldots, a_{m,0}^-,a_{1,1}^+, \ldots , a_{p,1}^+, a_{1,1}^-, \ldots, a_{q,1}^-)\neq 0.
$$
From the linearity of $f$ we can assume that  $a_{i,0}^+,a_{i,0}^- , a_{i,1}^+,a_{i,1}^-\in A_1 \cup \cdots \cup A_s\cup J$.
Since $A_iA_j=0$ for $i \neq j$, from the property of the $\ast$-graded exponent we have
$$
a_{1,0}^+, \ldots , a_{n,0}^+,a_{1,0}^-, \ldots, a_{m,0}^-,a_{1,1}^+, \ldots , a_{p,1}^+, a_{1,1}^-, \ldots, a_{q,1}^- \in A_{j_1}\oplus \cdots \oplus A_{j_k} + J
$$
for some $A_{j_1}, \ldots , A_{j_k}$ such that  $\mathrm{dim}_F(A_{j_1}\oplus \cdots \oplus A_{j_k})\leq d$ .
Thus $f$ is not an identity for one of the algebras $B_1, \ldots , B_t, D$. Hence $f \not\in Id_{\mathbb{Z}_2}^\ast(B_1 \oplus \cdots \oplus B_t\oplus D)$. In conclusion
$$
\mathrm{var}_{\mathbb{Z}_2}^\ast(A)\subseteq \mathrm{var}_{\mathbb{Z}_2}^\ast (B_1\oplus \cdots \oplus B_t \oplus D)
$$
and the proof is complete.

\bigskip

An application of Theorem \ref{decomposizione} is given in
 terms of $\ast$-graded codimensions.

\bigskip

\begin{cor} \label{supercodimensioni}
Let $A$ be a finitely generated PI-$\ast$-superalgebra. Then there exists a finite number of reduced $\ast$-superalgebras $B_1,
\ldots, B_t$ such that
$$
c_n^{(\mathbb{Z}_2,\ast)}(A)\simeq c_n^{(\mathbb{Z}_2,\ast)} (B_1\oplus \cdots \oplus B_t)
$$
\end{cor}
\noindent {\bf Proof.} By Theorem \ref{decomposizione} there is a finite number of reduced $\ast$-superalgebras $B_1, \ldots, B_t$ and a finite dimensional
$\ast$-superalgebra $D$ such that
$$
\mathrm{var}_{\mathbb{Z}_2}^\ast(A)=\mathrm{ var}_{\mathbb{Z}_2}^\ast (B_1\oplus \cdots \oplus B_t \oplus D)
$$

\noindent with
$\mathrm{exp}_{\mathbb{Z}_2}^\ast(A)=\mathrm{exp}_{\mathbb{Z}_2}^\ast(B_1)=\cdots=\mathrm{exp}_{\mathbb{Z}_2}^\ast(B_t)$
 and
$\mathrm{exp}_{\mathbb{Z}_2}^\ast(D)<\mathrm{exp}_{\mathbb{Z}_2}^\ast(A)$.
By Lemma \ref{codimensioni},
$$
c_n^{(\mathbb{Z}_2,\ast)}(B_1\oplus \cdots \oplus B_t) \leq c_n^{(\mathbb{Z}_2,\ast)}(B_1\oplus \cdots \oplus B_t \oplus D) \leq c_n^{(\mathbb{Z}_2,\ast)}(B_1\oplus \cdots
\oplus B_t)+c_n^{(\mathbb{Z}_2,\ast)}(D).
$$

\noindent Recalling that
$\mathrm{exp}_{\mathbb{Z}_2}^\ast(D)<\mathrm{exp}_{\mathbb{Z}_2}^\ast(B_1)=\mathrm{exp}_{\mathbb{Z}_2}^\ast(B_1\oplus
\cdots \oplus B_t)$ we have that
$$
c_n^{(\mathbb{Z}_2,\ast)}(A)\simeq c_n^{(\mathbb{Z}_2,\ast)}(B_1\oplus \cdots \oplus B_t)
$$
and the proof of the corollary is complete.

\bigskip

The following results give us a characterization of the varieties of
 $\ast$-superalgebras satisfying a Capelli identity. Let's start with the following lemma

\bigskip

\begin{lemma}  \label{capelliidentity}
Let $M^+$, $M^-$, $L^+$ and $L^-$ be natural numbers. If $A$ is a $\ast$-superalgebra satisfying
 the $\ast$-graded Capelli polynomials $Cap^{(\mathbb{Z}_2,\ast)}_{M^+}
[Y^+,X]$, $Cap^{(\mathbb{Z}_2,\ast)}_{M^-}
[Y^-,X]$, $Cap^{(\mathbb{Z}_2,\ast)}_{L^+} [Z^+,X]$ and $Cap^{(\mathbb{Z}_2,\ast)}_{L^-} [Z^-,X]$, then $A$ satisfies the Capelli
identity  $Cap_k (x_1, \ldots , x_{k} ; \bar{x}_1, \ldots ,
\bar{x}_{k-1})$, where $k=M^++M^-+L^++L^-$.
\end{lemma}
\noindent {\bf Proof.} Let $k=M^++M^-+L^++L^-$, then we obtain immediately the thesis if we
observe that
$$
Cap_{k} (x_1, \ldots , x_{k} ; \bar{x}_1, \ldots , \bar{x}_{k-1})=
$$
$$
Cap_{k} (\frac{y_{1}^++y_{1}^-}{2}+\frac{z_{1}^++z_{1}^-}{2},
\ldots ,
\frac{y_{k}^++y_{k}^-}{2}+\frac{z_{k}^++z_{k}^-}{2} ;
\bar{x}_1, \ldots , \bar{x}_{k-1})
$$

\noindent is a linear combinations of $\ast$-graded Capelli polynomials
alternating or in $m^+\geq M^+$ symmetric variables of zero degree, or in $m^-\geq M^-$ skew variables of zero degree, or in $l^+\geq L^+$
symmetric variables of one degree or in $l^-\geq L^-$ skew variables of one degree.

\bigskip

\begin{theorem} \label{finitelygenerated}
Let $\mathcal{V}_{\mathbb{Z}_2}^\ast$ be a variety of $\ast$-superalgebras. If $\mathcal{V}_{\mathbb{Z}_2}^\ast$
satisfies the Capelli identity of some rank, then
$\mathcal{V}_{\mathbb{Z}_2}^\ast=\mathrm{var}_{\mathbb{Z}_2}^\ast(A)$, for some finitely generated
$\ast$-superalgebra $A$.
\end{theorem}
\noindent {\bf Proof.} The proof follows very closely the proof given in \cite[Theorem 11.4.3]{GZ} for superalgebras.

\bigskip

\section{The $\ast$-superalgebra $UT_{\mathbb{Z}_2}^\ast(A_1, \ldots, A_m)$ }

\bigskip

In this section we recall the construction of the
$\ast$-superalgebra $UT_{\mathbb{Z}_2}^\ast(A_1, \ldots, A_m)$ given
in  section 3 of \cite{DSS} and we investigate the link between the
degrees of the $\ast$-graded Capelli polynomials and  the
$\ast$-graded identities of this $\ast$-superalgebra.

If $F$ is an algebraically closed field of characteristic zero, then, up to graded isomorphisms, the only finite dimensional simple $\ast$-superalgebras are the
following (see \cite[Theorem 7.6]{GSV})

\begin{enumerate}
  \item[$(1)$] $(M_{h,l}, \diamond)$, with $h\geq l \geq 0$, $h\neq 0$;
  \item[$(2)$] $(M_{h,l}\oplus M_{h,l}^{op}, exc)$, with $h\geq l \geq 0$, $h\neq 0$, and induced grading;
  \item[$(3)$] $(M_n+cM_n, \star)$, with involution given by $(a+cb)^\star=a^\diamond-cb^\diamond$;
  \item[$(4)$] $(M_n+cM_n, \dagger)$, with involution given by $(a+cb)^\dagger=a^\diamond+cb^\diamond$;
  \item[$(5)$] $((M_n+cM_n)\oplus (M_n+cM_n)^{op}, exc)$, with grading $(M_n\oplus M_n^{op},c(M_n\oplus M_n^{op}) )$;
\end{enumerate}

\noindent
where $\diamond =t,s$ denotes the transpose or symplectic involution and $exc$ is the exchange involution.
 Remember that the symplectic involution can occur only when $h=l$. Moreover $M_h=M_h(F)$ is the superalgebra of $h\times h$ matrices over $F$ with trivial grading,
  $M_{h,l}=M_{h+l}(F)$ is the superalgebra with grading
$\left ( \left ( \begin{array}{cc}
                F_{11} & 0 \\
                 0 & F_{22} \\
                  \end{array}
\right),
         \left ( \begin{array}{cc}
                    0 & F_{12} \\
                 F_{21} & 0 \\
\end{array} \right ) \right )$, where $F_{11}$, $F_{12}$, $F_{21}$, $F_{22}$ are $h \times h$, $h \times l$,
$l \times h$ and $l \times l$ matrices respectively, $h\geq l \geq 0$, $h\neq 0$ and $M_n+cM_n=M_n(F\oplus cF)$ denotes the simple superalgebra with grading $(M_n(F),
cM_n(F))$, where $c^2=1$.

\smallskip
Let $(A_1, \ldots, A_m)$ be a $m$-tuple of finite dimensional simple
$\ast$-superalgebras. For every $k =1, \ldots, m,$ the size of $A_k$
is given by
$$
s_k=\left\{\begin{array}{ll}
               h_k+l_k & \mathrm{if}\,\,  A_k=M_{h_k,l_k} \,\, \mathrm{or} \,\, A_k=M_{h_k,l_k}\oplus M_{h_k,l_k}^{op}; \\
               2n_k & \mathrm{if}\,\,  A_k=M_{n_k}+c M_{n_k} \,\, \mathrm{or} \,\, A_k=(M_{n_k}+c M_{n_k})\oplus (M_{n_k}+c M_{n_k})^{op} \\
             \end{array}
\right.
$$
and, set $\eta_0=0$, let $\eta_k=\Sigma_{i=1}^ks_i$ and
$\mathrm{B}l_k=\{\eta_{k-1}+1, \ldots , \eta_k\}$. Moreover, we denote
by $\gamma_m$ the orthogonal involution defined on the matrix
algebra $M_m$ by sending each $a \in M_m$ into the element
$a^{\gamma_m}\in M_m$ obtained reflecting $a$ along its secondary
diagonal. In particular for any matrix unit $e_{i,j}$ of $M_m$,
$e_{i,j}^{\gamma_m}=e_{m-j+1,m-i+1}$.

Then, we have a monomorphism of $\ast$-algebra
$$
\Delta : \bigoplus_{k=1}^m A_k \rightarrow
(M_{2\eta_m}, \gamma_{2\eta_m})
$$

defined by

$$
(c_1, \ldots, c_m)\rightarrow \left(
                                \begin{array}{cccccc}
                                  \bar{a}_1 &  &  &  &  & \\
                                  & \ddots &  &  &  &  \\
                                  &  &  \bar{a}_m &  &  &  \\
                                  &  &  & \bar{b}_m &  &  \\
                                  &  &  &  & \ddots &  \\
                                  &  &  &  &  & \bar{b}_1 \\
                                \end{array}
                              \right)
$$

\medskip
\noindent
where the elements $\bar{a}_i$ and  $\bar{b}_i$ are defined as follows:
\medskip

\noindent $\bullet$ if $c_i \in (M_{h,l}; \diamond)$, then
$\bar{a}_i=c_i$ and $\bar{b}_i=(c_i^{\diamond})^{\gamma_{h+l}}$;

\noindent $\bullet$ if $c_i=(a_i,b_i) \in (M_{h,l}\oplus M_{h,l}^{op}, exc)$, then $\bar{a}_i=a_i$ and $\bar{b}_i=b_i^{\gamma_{h+l}}$;

\noindent $\bullet$ if $c_i=a_i+cb_i \in (M_{n}+c M_{n}, \star)$, then $\bar{a}_i=\left(
                                                                                    \begin{array}{cc}
                                                                                      a_i & b_i \\
                                                                                      b_i & a_i \\
                                                                                    \end{array}
                                                                                  \right)$ and $\bar{b}_i=(\bar{a}_i^\bot)^{\gamma_{2n}}$ where $\left(
                                                                                    \begin{array}{cc}
                                                                                      x & y \\
                                                                                      y & x \\
                                                                                    \end{array}
                                                                                  \right)^\bot$ $=\left(
                                                                                    \begin{array}{cc}
                                                                                      x^{\diamond} & -y^{\diamond} \\
                                                                                      -y^{\diamond} & x^{\diamond} \\
                                                                                    \end{array}
                                                                                  \right)$;

\noindent $\bullet$ if $c_i=a_i+cb_i \in (M_{n}+c M_{n}, \dagger)$, then $\bar{a}_i=\left(
                                                                                    \begin{array}{cc}
                                                                                      a_i & b_i \\
                                                                                      b_i & a_i \\
                                                                                    \end{array}
                                                                                  \right)$ and $\bar{b}_i=(\bar{a}_i^\top)^{\gamma_{2n}}$ where $\left(
                                                                                    \begin{array}{cc}
                                                                                      x & y \\
                                                                                      y & x \\
                                                                                    \end{array}
                                                                                  \right)^\top$ $=\left(
                                                                                    \begin{array}{cc}
                                                                                      x^{\diamond} & y^{\diamond} \\
                                                                                      y^{\diamond} & x^{\diamond} \\
                                                                                    \end{array}
                                                                                  \right)$;

\noindent $\bullet$ if $c_i=(a_i+cb_i, u_i+cv_i) \in ((M_{n}+c M_{n})\oplus (M_{n}+c M_{n})^{op}, exc)$, then $\bar{a}_i=\left(
                                                                                    \begin{array}{cc}
                                                                                      a_i & b_i \\
                                                                                      b_i & a_i \\
                                                                                    \end{array}
                                                                                  \right)$ and $\bar{b}_i=\left(
                                                                                    \begin{array}{cc}
                                                                                      u_i & v_i \\
                                                                                      v_i & u_i \\
                                                                                    \end{array}
                                                                                  \right)^{\gamma_{2n}}$.

\medskip
Let denote by $D\subseteq (M_{2\eta_m},\gamma_{2\eta_m})$ the
$\ast$-algebra image of $\bigoplus_{i=1}^m A_i$ by $\Delta$ and set

\medskip
$$
V=\left(
  \begin{array}{cccccccc}
    0 & V_{12} & \cdots & V_{1m} &  &  &  &  \\
      & \ddots & \ddots & \vdots &  &  &  &  \\
      &   & 0 & V_{m-1m} &  &  &  &  \\
      &   &   & 0 &  &  &   &   \\
      &   &   &   & 0 & V_{mm-1} & \cdots & V_{m1} \\
      &   &   &   &   & \ddots & \ddots & \vdots \\
      &   &   &   &   &   & 0 & V_{21} \\
      &   &   &   &   &   &   & 0 \\
  \end{array}
\right)\subseteq M_{2\eta_m}
$$

\medskip
\noindent where, for $1 \leq i,j \leq m$, $i\neq j$, $V_{ij}=M_{s_i \times s_j}=M_{s_i \times s_j}(F)$ is the algebra of $s_i \times s_j$ matrices of $F$. Let define

$$
UT^\ast (A_1, \ldots, A_m)=D\oplus V \subseteq M_{2\eta_m}.
$$

\medskip
It is easy to see that $UT^\ast (A_1, \ldots, A_m)$ is a subalgebra
with involution of $(M_{2\eta_m}(F),\gamma_{2\eta_m})$ whose
Jacobson radical coincides with $V$.

Now, for any $m$-tuple $\tilde{g}=(g_1, \ldots , g_m)\in \mathbb{Z}_2^m$, we consider the map
$$
\alpha_{\tilde{g}} : \{1, \ldots ,2\eta_m\}\rightarrow \mathbb{Z}_2, \,\,\, i \rightarrow \left\{
                                                                               \begin{array}{ll}
                                                                                 \alpha_k (i-\eta_{k-1})+g_k & 1\leq i \leq \eta_m; \\
                                                                                 \alpha_k (2\eta_m-i+1-\eta_{k-1})+g_k & \eta_{m}+1\leq i \leq2 \eta_{m}.
                                                                               \end{array}
                                                                             \right.
$$
where $k \in \{ 1, \ldots , m\}$ is the (unique) integer such that
$i \in \mathrm{B}l_k$ and $\alpha_k$'s are maps  so defined:

\noindent $\cdot$  if $A_k \simeq M_{h,l}$ or $A_k  \simeq M_{h,l}\oplus M_{h,l}$, then
$$
\alpha_k : \{ 1, \ldots , h+l\} \rightarrow \mathbb{Z}_2,  \,\,\, \alpha_k(i)= \left\{
                                                                               \begin{array}{ll}
                                                                                 0 & 1\leq i \leq h; \\
                                                                                1 & h+1\leq i \leq h+l.
                                                                               \end{array}
                                                                             \right.
$$

\noindent $\cdot$  if $A_k \simeq M_{n}+cM_n$ or $A_k  \simeq (M_{n}+cM_n)\oplus (M_{n}+cM_n)$, then
$$
\alpha_k : \{ 1, \ldots ,2n\} \rightarrow \mathbb{Z}_2,  \,\,\, \alpha_k(i)= \left\{
                                                                               \begin{array}{ll}
                                                                                 0 & 1\leq i \leq n; \\
                                                                                1 & n+1\leq i \leq 2n.
                                                                               \end{array}
                                                                             \right.
$$

The map $\alpha_{\tilde{g}}$ induces an elementary grading on
$UT^\ast (A_1, \ldots, A_m)$ with respect to which
$\gamma_{2\eta_m}$ is a graded involution. We shall use the symbol

$$
UT_{\mathbb{Z}_2, \tilde{g}}^\ast (A_1, \ldots, A_m)
$$

\medskip
\noindent to indicate
the $\ast$-superalgebra defined by the $m$-tuple $\tilde{g}$. We
observe that the $k$-th simple component of the maximal semisimple
$\ast$-graded subalgebra of this $\ast$-superalgebra is isomorphic
to $A_k$. When convenient, any such $\ast$-superalgebra is simply
denoted by

$$
UT_{\mathbb{Z}_2}^\ast (A_1, \ldots, A_m).
$$

\bigskip

In the next lemma we establish the link between the degrees of the
$\ast$-graded Capelli polynomials and the $\ast$-graded polynomial
identities of $UT_{\mathbb{Z}_2, \tilde{g}}^\ast (A_1, \ldots,
A_m)$.
 For all $i=1, \ldots , m$, we write $$ A_i= A_{i,0}^+ \oplus
A_{i,0}^- \oplus A_{i,1}^+ \oplus A_{i,1}^- .$$ Let
$(d_{0}^{\pm})_i=\mathrm{dim}_FA_{i,0}^{\pm}$ and
$(d_{1}^{\pm})_i=\mathrm{dim}_FA_{i,1}^{\pm},$ if we set
$d_0^{\pm}:=\sum_{i=1}^m (d_{0}^{\pm})_i $ and
$d_1^{\pm}:=\sum_{i=1}^m (d_{1}^{\pm})_i, $ then we have the
following

\bigskip

\begin{lemma}\label{capelliidentita} Let $\tilde{g}=(g_1, \ldots , g_{m})$ be a fixed element of $\mathbb{Z}_2^m$ and $A=UT_{\mathbb{Z}_2, \tilde{g}}^\ast (A_1, \ldots,
A_m)$, with  $A_i$ finite dimensional simple $\ast$-superalgebra.  Let $0< \bar{m} \leq m$ denote the number of the finite dimensional simple $\ast$-superalgebras with trivial
grading.

\begin{enumerate}
  \item [1.] If $\bar{m}=0$, $Cap_{q^+}^{(\mathbb{Z}_2, \ast)} [Y^+,X],\,
Cap_{q^-}^{(\mathbb{Z}_2, \ast)} [Y^-,X], \,
Cap_{k^+}^{(\mathbb{Z}_2, \ast)} [Z^+,X] \,  \textrm{and}  \,
 Cap_{k^-}^{(\mathbb{Z}_2, \ast)} [Z^-,X]$
 are in
$Id^{\ast}_{\mathbb{Z}_2}(A)$ if and only if
$q^+\geq d_0^++m$, $q^-\geq d_0^- +m$, $k^+\geq d_1^++m$ and $k^-\geq d_1^-+m$;
  \item [2.]  If $0<\bar{m}\leq m$, let  $\tilde{m}$  be the number of blocks of consecutive $\ast$-superalgebras with trivial grading that appear in $(A_1, \ldots,
A_m)$. Then
 $Cap_{q^+}^{(\mathbb{Z}_2, \ast)} [Y^+,X],
Cap_{q^-}^{(\mathbb{Z}_2, \ast)} [Y^-,X],$  $
Cap_{k^+}^{(\mathbb{Z}_2, \ast)} [Z^+,X]$   and $Cap_{k^-}^{(\mathbb{Z}_2, \ast)} [Z^-,X]$
 are in $Id^{\ast}_{\mathbb{Z}_2}(A)$ if and only if $q^+> d_0^++(m-\bar{m})+(\tilde{m}-1)+r_0$, $q^->d_0^- +(m-\bar{m})+(\tilde{m}-1)+r_0$, $k^+> d_1^++(m-\bar{m})+(\tilde{m}-1)+r_1$ and $k^-> d_1^-+(m-\bar{m})+(\tilde{m}-1)+r_1$,
   where $r_0$, $r_1$ are two non negative integers depending on the grading $\tilde{g}$, with $r_0+r_1=\bar{m}-\tilde{m}$.
\end{enumerate}

\end{lemma}
\noindent {\bf Proof.}  We will prove the statement only for
$Cap_{q^+}^{(\mathbb{Z}_2, \ast)} [Y^+,X]$ the $\ast$-graded Capelli
polynomial  alternating on $q^+$ symmetric variables of degree zero
since on the other cases the proofs are similar.

$1.$ Let $\bar{m}=0$. To prove the necessary condition of the statement for the symmetric variables of degree zero it is sufficient to prove that
$Cap^{(\mathbb{Z}_2, \ast)}_{q^+}[Y^+,X]$ is not in $Id^{\ast}_{\mathbb{Z}_2}(A)$ when $q^+=d_0^++m-1$.

We start considering separately the components $A_i$ of $A$. In each
$\ast$-superalgebra $A_i$ we can take $(d_{0}^+)_i$ symmetric
elements of homogeneous degree zero
$$
S_i=\{s_{\alpha_{i-1}+ i}, \ldots , s_{\alpha_i+i-1}\}
$$
for $i=1, \ldots , m$, where $\alpha_0=0$ and $\alpha_i=\sum_{j=0}^i(d_0^+)_j$ and a set of elements of $A_i$
$$
U_i=\{a_{\alpha_{i-1}+ i}, \ldots , a_{\alpha_i+i-2}\}
$$
such that
$$
Cap^{(\mathbb{Z}_2, \ast)}_{(d_0^+)_i}(s_{\alpha_{i-1}+ i}, \ldots , s_{\alpha_i+i-1};a_{\alpha_{i-1}+ i}, \ldots , a_{\alpha_i+i-2})
=$$
$$
\left\{
        \begin{array}{ll}
            e_{r_i,s_i} & \hbox{if} \,\,\, (M_{h_i,l_i}, \diamond);\\
            (e_{r_i,s_i},0) & \hbox{if} \,\,\,(M_{h_i,l_i}\oplus M_{h_i,l_i}^{op}, exc);\\
            e_{r_i,s_i} & \hbox{if} \,\,\, (M_{n_i}+cM_{n_i}, \star)\,\,\, \hbox{or} \,\,\, (M_{n_i}+cM_{n_i}, \dagger);\\
            ((e_{r_i,s_i},0),(0,0)) & \hbox{if} \,\,\, ((M_{n_i}+cM_{n_i})\oplus (M_{n_i}+cM_{n_i})^{op}, exc),
        \end{array}
\right.
$$

\medskip\noindent
where $\diamond =t,s$ denotes the transpose or symplectic involution, $exc$ is the exchange involution, $(a+cb)^\star=a^\diamond-cb^\diamond$ and
$(a+cb)^\dagger=a^\diamond+cb^\diamond$.

For any $1\leq i \leq m$, if $\phi_i$ is the $\ast$-embedding of $A_i$ in A, then let

$$
\bar{S}_i=\{\bar{s}_{\alpha_{i-1}+ i}, \ldots , \bar{s}_{\alpha_i+i-1}\}
$$
and
$$
\bar{U}_i=\{\bar{a}_{\alpha_{i-1}+ i}, \ldots , \bar{a}_{\alpha_i+i-2}\}
$$

\noindent
denote the images of $S_i$ and $U_i$ by $\phi_i$, respectively.

Let observe that in $A$ we can consider appropriate symmetric elements of homogeneous degree zero in $J_0^+$
$$
\bar{s}_{\alpha_i+i}=e_{h,k}+e_{h,k}^\ast
$$
and elementary matrices of $A$
$$
\bar{a}_{\alpha_{i}+ i-1}=e_{s_i,h} \,\,\, \mathrm{and} \,\,\, \bar{a}_{\alpha_{i}+ i}=e_{k,r_{i+1}}
$$
such that
$$
Cap^{(\mathbb{Z}_2, \ast)}_{(d_{0}^+)_i}(\bar{s}_{\alpha_{i-1}+i}, \ldots , \bar{s}_{\alpha_{i}+i-1};  \bar{a}_{\alpha_{i-1}+i},\ldots  ,\bar{a}_{\alpha_{i}+i-2})
\bar{a}_{\alpha_{i}+i-1}\bar{s}_{\alpha_{i}+i}\bar{a}_{\alpha_{i}+i}
$$
$$
 Cap^{(\mathbb{Z}_2, \ast)}_{(d_{0}^+)_{i+1}}(\bar{s}_{\alpha_{i}+(i+1)}, \ldots , \bar{s}_{\alpha_{i+1}+i};  \bar{a}_{\alpha_{i}+(i+1)},\ldots ,
 \bar{a}_{\alpha_{i+1}+(i-1)})\neq
0.
 $$
From now on, we will put $Cap^{(\mathbb{Z}_2, \ast)}_{(d_{0}^+)_i}=Cap^{(\mathbb{Z}_2, \ast)}_{(d_{0}^+)_i}(\bar{s}_{\alpha_{i-1}+i}, \ldots , \bar{s}_{\alpha_{i}+i-1};
\bar{a}_{\alpha_{i-1}+i},\ldots  ,\bar{a}_{\alpha_{i}+i-2})$. It follows that
$$
Cap^{(\mathbb{Z}_2, \ast)}_{q^+}(\bar{s}_{1}, \ldots ,  \bar{s}_{\alpha_m+(m-1)};
 \bar{a}_{1},\ldots ,\bar{a}_{\alpha_m+(m-2)})=
 $$
 $$
Cap^{(\mathbb{Z}_2, \ast)}_{(d_{0}^+)_{1}}\bar{a}_{\alpha_1}\bar{s}_{\alpha_1+1}\bar{a}_{\alpha_1+1}
Cap^{(\mathbb{Z}_2, \ast)}_{(d_{0}^+)_{2}}
\cdots \cdots
Cap^{(\mathbb{Z}_2, \ast)}_{(d_{0}^+)_{m-1}}\bar{a}_{\alpha_{m-1}}\bar{s}_{\alpha_{m-1}+1}\bar{a}_{\alpha_{m-1}+1}Cap^{(\mathbb{Z}_2, \ast)}_{(d_{0}^+)_{m}}\neq 0.
$$

Conversely, let $q^+\geq  d_0^++m$. We observe that any monomial of
elements of $A$ containing at least $m$ elements of $J_0^+$ must be
zero. Then we claim that any multilinear polynomial $
\tilde{f}=\tilde{f}(y_1, \ldots , y_{d_0^++m};x_1, x_2, \ldots ) $
alternating  on $d_0^++m$  symmetric variables of degree zero must
vanish in $A$. In fact, by multilinearity, we can consider only
substitutions $\varphi: y_i^+ \rightarrow \bar{s}_i$, $x_i
\rightarrow \bar{a}_i$ such that $\bar{s}_i \in D_0^+ \cup J_0^+$
for $1 \leq i \leq d_0^++m$.

However, since $\mathrm{dim}_FD_{0}^+=d_0^+$, if we substitute at least $d_0^++1$ variables in elements of $D_0^+ $ the polynomial vanishes. On the other hands, if we
substitute at least $m$ elements of $J_0^+$, we also get that $\tilde{f}$ vanishes in $A$. The outcome  of this is that $A$ satisfies  $Cap_{d_0^++m}^{(\mathbb{Z}_2, \ast)}
[Y^+,X]$ and so $Cap_{q^+}^{(\mathbb{Z}_2, \ast)} [Y^+,X]$, with $q^+\geq d_0^++m$.

\medskip
$2.$ First let assume that $\bar{m}=m$. We recall that
$$
UT^\ast (A_1, \ldots, A_m)=D\oplus V \subseteq M_{2\eta_m},
$$
where $D\subseteq (M_{2\eta_m},\gamma_{2\eta_m})$ the $\ast$-algebra image of $\bigoplus_{i=1}^m A_i$ by $\Delta$ and

$$
V=\left(
  \begin{array}{cccccccc}
    0 & V_{12} & \cdots & V_{1m} &  &  &  &  \\
      & \ddots & \ddots & \vdots &  &  &  &  \\
      &   & 0 & V_{m-1m} &  &  &  &  \\
      &   &   & 0 &  &  &   &   \\
      &   &   &   & 0 & V_{mm-1} & \cdots & V_{m1} \\
      &   &   &   &   & \ddots & \ddots & \vdots \\
      &   &   &   &   &   & 0 & V_{21} \\
      &   &   &   &   &   &   & 0 \\
  \end{array}
\right)\subseteq M_{2\eta_m}
$$

\medskip Notice that, for a fixed $\tilde{g}=(g_1, \ldots , g_{m})\in
\mathbb{Z}_2^m,$ if $g_i = g_j$, $1\leq i,j \leq m$, then the
elements of the blocks $V_{i,j}$ are homogeneous of degree zero,
otherwise, if $g_i\neq g_j$, they are homogeneous of degree one.
Suppose that in  $\tilde{g}=(g_1, \ldots , g_{m})$ there are $p \ge
1$ different string of zero and one,  i.e.
$$
\tilde{g}=(g_1, \ldots , g_{t_1}, g_{t_1+1}, \ldots , g_{t_1+t_2},\ldots , g_{t_1+\cdots + t_{p-1}+1}, \ldots , g_{t_1+\cdots +t_p}),
$$
where $t_1+\cdots +t_p=m$,
$$
g_1=\cdots =g_{t_1},$$
$$g_{t_1+1}=\cdots =g_{t_1+t_2},$$
$$\ldots \ldots$$
$$g_{t_1+\cdots + t_{p-1}+1}=\cdots =g_{t_1+\cdots +t_p}
$$
and
$$g_{t_1+\cdots + t_{i}}\neq g_{t_1+\cdots + t_{i}+1},$$
$\forall \,\, i=1,\ldots , p-1$.

 As in the previous case  we can find in  $A$  symmetric elements of degree zero

$$
\bar{S}_i=\{\bar{s}_{\alpha_{i-1}+ i}, \ldots , \bar{s}_{\alpha_i+i-1}, \bar{s}_{\alpha_i+i}\}
$$
and generic elements
$$
\bar{U}_i=\{\bar{a}_{\alpha_{i-1}+ i}, \ldots , \bar{a}_{\alpha_i+i-2},\bar{a}_{\alpha_i+i-1}, \bar{a}_{\alpha_i+i}\}
$$
such that, $\forall \,\, i=1,\ldots p$,

$$
Cap^{(\mathbb{Z}_2, \ast)}_{q_{i}}(\bar{s}_{\alpha_{\tilde{t}_{i-1}}+ (\tilde{t}_{i-1}+1)}, \ldots , \bar{s}_{\alpha_{\tilde{t}_{i}}+(\tilde{t}_{i}-1)};
\bar{a}_{\alpha_{\tilde{t}_{i-1}}+ (\tilde{t}_{i-1}+1)}, \ldots , \bar{a}_{\alpha_{\tilde{t}_{i}}+(\tilde{t}_{i}-2)})=
$$
$$
Cap^{(\mathbb{Z}_2, \ast)}_{(d_{0}^+)_{\tilde{t}_{i-1}+1}}
\bar{a}_{\alpha_{(\tilde{t}_{i-1}+1)}+ \tilde{t}_{i-1}}\bar{s}_{\alpha_{(\tilde{t}_{i-1}+1)}+ (\tilde{t}_{i-1}+1)}\bar{a}_{\alpha_{
(\tilde{t}_{i-1}+1)}+ (\tilde{t}_{i-1}+1)}
Cap^{(\mathbb{Z}_2, \ast)}_{(d_{0}^+)_{\tilde{t}_{i-1}+2}}$$
$$
\cdots \cdots \cdots
\cdots  Cap^{(\mathbb{Z}_2, \ast)}_{(d_{0}^+)_{\tilde{t}_{i}}}=b_i\neq 0,
$$

\noindent where $\tilde{t}_0=t_0=0$, $\tilde{t}_i=\sum_{j=0}^it_j$ and $q_{i}=(d_{0}^+)_{\tilde{t}_{i-1}+1}+\cdots +(d_{0}^+)_{\tilde{t}_{i}}+(t_{i}-1)$.

\medskip
Furthermore we can find in $A$ elementary  matrices $ E_1, \ldots ,
E_{p-1}, $ such that

$$
Cap^{(\mathbb{Z}_2, \ast)}_{d_{0}^++m-p}=
Cap^{(\mathbb{Z}_2, \ast)}_{q_{1}}E_1
Cap^{(\mathbb{Z}_2, \ast)}_{q_{2}}E_2 \cdots  Cap^{(\mathbb{Z}_2, \ast)}_{q_{p-1}} E_{p-1} Cap^{(\mathbb{Z}_2, \ast)}_{q_{p}}=
$$
$$
b_1E_1b_2E_2\cdots b_{p-1} E_{p-1}b_p\neq 0.
$$

This implies that, for $r_0=m-p$,

$$
Cap_{d_0^++r_0}^{(\mathbb{Z}_2, \ast)}[Y^+,X]\notin Id^{\ast}_{\mathbb{Z}_2}(A).
$$
Moreover, let's observe that any monomial of elements of $A$
containing at least $r_0+1=(m-p)+1$ elements of $J_0$ must be zero.
Then, similarly to the previous case, we obtain that $A$ satisfies
$Cap_{d_0^++r_0+1}^{(\mathbb{Z}_2, \ast)}[Y^+,X]$.

If $0<\bar{m}<m$, let $\tilde{m}$  be the number of blocks of consecutive $\ast$-superalgebras with trivial grading that appear in $(A_1, \ldots,
A_m)$.   By considering separately the blocks of consecutive $\ast$-superalgebras with trivial and non-trivial grading and by using arguments  similar to those of the proof of case 1, it easily follows that
 $Cap_{q^+}^{(\mathbb{Z}_2, \ast)} [Y^+,X],\,
Cap_{q^-}^{(\mathbb{Z}_2, \ast)} [Y^-,X], \,
Cap_{k^+}^{(\mathbb{Z}_2, \ast)} [Z^+,X]$   and $Cap_{k^-}^{(\mathbb{Z}_2, \ast)} [Z^-,X]$
 are in $Id^{\ast}_{\mathbb{Z}_2}(A)$ if and only if $q^+> d_0^++(m-\bar{m})+(\tilde{m}-1)+r_0$, $q^->d_0^- +(m-\bar{m})+(\tilde{m}-1)+r_0$, $k^+> d_1^++(m-\bar{m})+(\tilde{m}-1)+r_1$ and $k^-> d_1^-+(m-\bar{m})+(\tilde{m}-1)+r_1$,
   where $r_0$, $r_1$ are two non negative integers depending on the grading $\tilde{g}$, with $r_0+r_1=\bar{m}-\tilde{m}$.

\bigskip

\section{Asymptotics for $\ast$-graded Capelli identities}

In this section we shall study
$\mathcal{U}=\mathrm{var}_{\mathbb{Z}_2}^\ast(\Gamma^\ast_{M^{\pm}+1,
L^{\pm}+1})$ and  we shall find a close relation among the
asymptotics of $c^{\ast}_n(\Gamma^{\ast}_{M^{\pm}+1,L^{\pm}+1})$ and
$c^{\ast}_n(A)$, where $A$  is a finite dimensional simple
$\ast$-superalgebra. Let
$$
R= {A} \oplus J
$$

\noindent where ${A}$ is a   finite dimensional simple  $\ast$-superalgebra
and $J=J(R)$ is its Jacobson radical.

From now on we put $M^{\pm} =\textrm{dim}_F {A_0}^{\pm}$ and
$L^{\pm}=\textrm{dim}_F {A_1}^{\pm}$.

Let's begin with some technical lemmas that hold for any finite
dimensional simple $\ast$-superalgebra ${A}$.

\begin{lemma} \label{radicaltrasposta} The Jacobson radical $J$ can be decomposed
into the direct sum of four ${A}$-bimodules
$$ J=J_{00}\oplus J_{01} \oplus J_{10} \oplus J_{11} $$ where, for
$p,q \in \{0,1\}$, $J_{pq}$ is a left faithful module or a $0$-left
module according to $p=1$, or $p=0$, respectively. Similarly,
$J_{pq}$ is a right faithful module or a $0$-right module according
to $q=1$ or $q=0$, respectively. Moreover, for $p,q,i,l \in
\{0,1\}$, $J_{pq}J_{ql}\subseteq J_{pl}$, $J_{pq}J_{il}=0$ for
$q\neq i$ and there exists a finite dimensional nilpotent
$\ast$-superalgebra $N$ such that $N$ commutes with ${A}$ and
$J_{11}\cong {A}\otimes_F N$ (isomorphism of ${A}$-bimodules and of
$\ast$-superalgebras).
\end{lemma}
\noindent {\bf Proof.} It follows  from  Lemma 2 in
\cite{GZ2} and  Lemmas 1,6 in \cite{BS}.

\bigskip
Notice that $J_{00}$ and $J_{11}$ are stable under the involution
whereas $J_{01}^\ast=J_{10}$.

\bigskip

\begin{lemma} \label{j10trasposta}
    If $\Gamma^\ast_{M^{\pm}+1,L^{\pm}+1} \subseteq Id^{\ast}_{\mathbb{Z}_2}(R)$, then $J_{10}=J_{01}=(0)$.
\end{lemma}

\noindent {\bf Proof.} By   Lemma \ref{capelliidentita}  we have
that ${A}$ does not satisfy  $Cap^{(\mathbb{Z}_2,
\ast)}_{M^+}[Y^+,X]$. Then there exist elements $a^+_1, \ldots ,
a^+_{M^+} \in {A_0}^{+}$ and $b_1, \ldots , b_{M^+-1} \in {A}$ such
that

$$
Cap^{(\mathbb{Z}_2, \ast)}_{M^+}(a^+_1, \ldots , a^+_{M^+} ;b_1, \ldots ,
b_{M^+-1})=
$$

$$
\left\{
\begin{array}{ll}
e_{1,h+l} & \hbox{if} \, {A}=(M_{h,l},\diamond),  \diamond=t,s; \\
\tilde{e}_{1,h+l} & \hbox{if} \, {A}=(M_{h,l}\oplus M_{h,l}^{op},exc); \\
e_{1,n} & \hbox{if}  \, {A}=(M_{n}+ c M_{n},\star) \, \hbox{or}   {A}=(M_{n}+ c M_{n},\dagger);\\
\tilde{e}_{1,n}  & \hbox{if} \, {A}=((M_{n}+ c M_{n}) \oplus ( M_{n}+
c M_{n})^{op},exc)
\end{array}
\right.
$$

\medskip\noindent where the $e_{i,j}$'s are the usual matrix units and
$\tilde{e}_{i,j}=(e_{i,j}, e_{j,i})$.  We write
$J_{10}=(J_{10})_0\oplus (J_{10})_1$ and $J_{01}=(J_{01})_0\oplus
(J_{01})_1$. Let $d_0 \in (J_{01})_0$, then $d_0^\ast \in
(J_{10})_0$ and $d_0+d_0^\ast \in (J_{01} \oplus J_{10})^{+}_0$.
Since $\Gamma^\ast_{M^{\pm}+1,L^{\pm}+1} \subseteq
Id^{\ast}_{\mathbb{Z}_2}(R)$ it follows that there exists $b_{M^{+}}\in {A}$ such that

$$
0=Cap^{(\mathbb{Z}_2, \ast)}_{M^{+}+1}(a^+_1, \ldots ,
a^+_{M^{+}},d_0+d_0^\ast;b_1, \ldots , b_{M^{+}-1}, b_{M^{+}})=
$$

$$ \left\{
\begin{array}{ll}
e_{1,h+l}d_0^\ast \pm  d_0e_{1,h+l}  & \hbox{if} \, {A}=(M_{h,l},\diamond), \, \diamond=t,s; \\
\tilde{e}_{1,h+l} d_0^\ast \pm d_0\tilde{e}_{1,h+l} & \hbox{if} \,  {A}=(M_{h,l}\oplus M_{h,l}^{op},exc); \\
{e}_{1,n} d_0^\ast \pm  d_0{e}_{1,n} & \hbox{if} \,   {A}=(M_{n}+ c M_{n},\star) \, \hbox{or} \,  {A}=(M_{n}+ c M_{n},\dagger); \\
\tilde{e}_{1,n} d_0^\ast \pm  d_0\tilde{e}_{1,n} & \hbox{if} \,
{A}=((M_{n}+ c M_{n})\oplus( M_{n}+ c M_{n})^{op},exc).
\end{array}
\right.
$$

\medskip \noindent  If ${A}=(M_{h,l},\diamond)$, then $ e_{1,h+l}d_0^\ast \pm d_0e_{1,h+l}=0$ and, so, $ e_{1,h+l}d_0^\ast
=\mp d_0e_{1,h+l} \in (J_{01})_0\cap (J_{10})_0=(0)$. Hence $d_0=0$, for all
$d_0 \in (J_{01})_0$. Thus $(J_{01})_0=(0)$ and $(J_{10})_0=(0)$.
Similarly for the other finite dimensional simple $\ast$-superalgebras we obtain that $(J_{01})_0=(J_{10})_0=(0)$.
Analogously it easy to show that $(J_{01})_1=(J_{10})_1=(0)$ and the lemma is proved.

\bigskip

\begin{lemma} \label{Ntrasposta}
    Let $J_{11} \cong {A} \otimes_F N$, as in Lemma
    \ref{radicaltrasposta}.   If
    $\Gamma^\ast_{M^{\pm}+1,L^{\pm}+1}\subseteq Id^{\ast}_{\mathbb{Z}_2}(R)$, then $N$ is commutative.
\end{lemma}
\noindent {\bf Proof.} Let $N$ be the finite dimensional nilpotent
$\ast$-superalgebra of  Lemma \ref{radicaltrasposta}. Write
$N=N_0^+\oplus N_0^- \oplus N_1^+ \oplus N_1^-$,  where $N_0^+,$
$N_0^-,$  $N_1^+$ and $N_1^-$ denote the subspaces of symmetric and
skew symmetric elements of $N$ of homogeneous degree 0 and $1$
respectively.

We shall prove that $N$ is commutative when $ A =
(M_{h,l},\diamond)$, with $\diamond= t$ or $s$.   Similar
calculations for the other finite dimensional simple
$\ast$-superalgebras lead to the same conclusion.

Let's start by proving that $N_0^\pm$ commutes with $N_i^{\pm}, i=
0,1$. Let  $e_1^+,\ldots , e_{M^+}^+$ be  a basis of ${A_0}^{+}$
with
$$
e_1^+ = \left\{
\begin{array}{ll}
e_{1,2} + e_{2,1}  & \hbox{if} \, {A}=(M_{h,l},t) \,; \\
e_{1,2} + e_{h+2,h+1}  & \hbox{if} \, {A}=(M_{h,h},s)
\end{array}
\right.
$$
and let $a_0=a_1=e_{2,1}$, $a_2, \ldots,  a_{M^+-1} \in {A}$ such
that $a_0e_1^+a_1e_2^+\cdots a_{M^+-1}e_{M^+}^+ = e_{2,h+l}$ and
 $a_0e_{\sigma(1)}^+a_1 \cdots a_{M^+-1}e_{\sigma(M^+)}^+ = 0$
 for any $\sigma \in S_{M^+}$, $\sigma \neq id$.
Let $ d_1 \in N_0^\pm$ and $e_0^+=(e_{1,2}\pm e_{1,2}^\diamond )
d_1$, with $\diamond= t$ or $s$. Since $N$ commutes with ${A}$ we
obtain that $e_0^+ \in {R}_0^+$. If we put
$\bar{a}_0=a_0d_2=e_{2,1}d_2$, with $d_2 \in N_i^\pm$, $i=0,1$, then
$$
0=Cap^{(\mathbb{Z}_2, \ast)}_{M^++1}(e_0^+,e_1^+,\ldots ,e_M^+;\bar{a}_0, a_1, \ldots , a_{M^+-1} )= [d_1, d_2]e_{1,h+l}
$$

\noindent and so $[d_1,d_2]= 0$ for all $d_1 \in N_0^\pm,  d_2 \in
N_i^{\pm}, \, i =0,1.$

Let's now prove that $N_1^\pm$ commutes with $N_1^{\pm}.$ Let
$e_1^+,\ldots , e_{M^+}^+$ be a basis of ${A_0}^{+},$ with
$$
e_1^+ = \left\{
\begin{array}{ll}
e_{1,1}   & \hbox{if} \, {A}=(M_{h,l},t) \,; \\
e_{1,1} + e_{h+1,h+1}  & \hbox{if} \, {A}=(M_{h,h},s)
\end{array}
\right.
$$
and let $a_0=e_{h+l,1}$, $a_1, a_2, \ldots,  a_{M^+-1} \in {A}$ such
that $a_0e_1^+a_1\cdots a_{M^+-1}e_{M^+}^+ = e_{h+l,1}$ (if
$\diamond$ =s then $h=l$) and
 $a_0e_{\sigma(1)}^+a_1 \cdots a_{M^+-1}e_{\sigma(M^+)}^+ = 0$
 for any $\sigma \in S_{M^+}$, $\sigma \neq id$.

\noindent Let $(e_{1, h+l}\pm {e}_{1,h+l}^\diamond) \in {A}_1^\pm$
and $d_1, d_2 \in N_1^\pm$ such that, for $i=1,2$, $c_i^+=(e_{1,
h+l}\pm {e}_{1,h+l}^\diamond) d_i$. Since $N$ commutes with ${A}$
then $c_i^+ \in {R}_0^+$, $i=1,2$. If $a_M=e_{1,1}$ then
$$
0=Cap^{(\mathbb{Z}_2, \ast)}_{M^++2}(c_1^+,e_1^+,\ldots ,e_M^+, c_2^+;\bar{a}_0, a_1, \ldots , a_{M^+-1}, a_M )= [d_1, d_2]e_{1,h+l}
$$
($h=l$ for $\diamond=s$) and so $[d_1,d_2] =0,$ for all $d_1,  d_2 \in N_1^{\pm}$ and we are done.

\bigskip

\begin{lemma} \label{esponente}

$
\mathrm{exp}_{\mathbb{Z}_2}^{\ast}(\mathcal{U})=M^++M^-+L^++L^-=M+L=\mathrm{exp}_{\mathbb{Z}_2}^{\ast}(A).
$
\end{lemma}
\noindent {\bf Proof.} By the definition of minimal variety (see
Definition 2.1 in \cite{DSS}) the $\ast$-graded exponent of
$\mathcal{U}$ is equal to the $\ast$-graded exponent of some minimal
variety of $\ast$-superalgebras lying in $\mathcal{U}$. Moreover, by
the classification of minimal varieties of PI-$\ast$-superalgebras
of finite basic rank given in \cite[Theorem 2.2]{DSS}, we have
$$
\mathrm{exp}_{\mathbb{Z}_2}^\ast(\mathcal{U})=
\mathrm{max}\{\mathrm{exp}_{\mathbb{Z}_2}^{\ast}(UT_{\mathbb{Z}_2}^\ast(A_1,
\dots ,A_m)) \, | \, UT_{\mathbb{Z}_2}^\ast(A_1, \dots
,A_m) \in \mathcal{U} \}.
$$
Then, by Lemma \ref{capelliidentita},
$$
\mathrm{exp}_{\mathbb{Z}_2}^\ast(\mathcal{U})\geq
\mathrm{exp}_{\mathbb{Z}_2}^\ast(UT_{\mathbb{Z}_2}^\ast ({A}))=M+L.
$$

\medskip
\noindent On the other hand, since
$\mathrm{exp}_{\mathbb{Z}_2}^{\ast}(UT_{\mathbb{Z}_2}^\ast(A_1,
\dots ,A_m)) =d_0^\pm + d_1^\pm$, we have that
$$\mathrm{exp}_{\mathbb{Z}_2}^\ast(\mathcal{U}) \leq M+L
$$

\noindent and the
proof is completed.

\bigskip

Now we are able to prove the main result.

\bigskip
\begin{theorem} \label{teorematrasposta}
For suitable natural numbers  $M^+$, $M^-$, $L^+$, $L^-$ there exists a finite dimensional simple $\ast$-superalgebra $A$ such that
$$
\mathcal{U}=\mathrm{var}_{\mathbb{Z}_2}^\ast(\Gamma^{\ast}_{M^{\pm}+1,
L^{\pm}+1})=\mathrm{var}_{\mathbb{Z}_2}^\ast(A\oplus D),
$$

\medskip
\noindent where $D$ is a finite dimensional $\ast$-superalgebra such
that $\mathrm{exp}_{\mathbb{Z}_2}^\ast(D)<M+L$, with $M=M^++M^-$ and
$L=L^++L^-$. In particular

\begin{enumerate}
  \item[1)] If $M^{\pm}=\frac{h(h\pm 1)}{2}+\frac{l(l\pm 1)}{2}$ and $L^{\pm}=hl$, with $h\geq l > 0$, then $A=(M_{h,l},t)$;
  \item[2)] If $M^{\pm}=h^2$ and $L^{\pm}=h(h\mp 1)$, with $h>0$, then $A=(M_{h,h},s)$;
  \item[3)] If $M^{\pm}=h^2+l^2$ and $L^{\pm}=2hl$, with $h\geq l > 0$, then $A=(M_{h,l}\oplus M_{h,l}^{op},exc)$;
  \item[4)] If $M^{+}=L^{\pm}=\frac{n(n+1)}{2}$, $M^{-}=L^{\mp}=\frac{n(n-1)}{2}$,  with $n>0$, then $A=(M_n+cM_n,\ast)$, where $(a+cb)^\ast= a^t \pm cb^t$;
  \item[5)]  If $M^{+}=L^{\pm}=\frac{n(n-1)}{2}$, $M^{-}=L^{\mp}=\frac{n(n+1)}{2}$,  with $n>0$, then $A=(M_n+cM_n,\ast)$, where $(a+cb)^\ast= a^s \pm cb^s$;
 \item [6)] If $M^{\pm}=L^{\pm}=n^2$,  with $n>0$,  then $A=((M_n+cM_n)\oplus (M_n+cM_n)^{op},exc)$.
\end{enumerate}

\end{theorem}

\noindent {\bf Proof.} By Lemma \ref{esponente} we have that
$\mathrm{exp}_{\mathbb{Z}_2}^{\ast}(\mathcal{U})=M+L$. Let $B$  be a
generating $\ast$-superalgebra of $\mathcal{U}$. From Theorem \ref{finitelygenerated} and by \cite{GIL2},
since any finitely generated $\ast$-superalgebra satisfies the same
$\ast$-graded polynomial identities of a finite-dimensional
$\ast$-superalgebra, we can assume
 that $B$ is  finite dimensional. Thus, by Theorem
\ref{decomposizione}, there exists a finite
 number of reduced $\ast$-superalgebras $B_1,\ldots , B_t$ and a finite dimensional $\ast$-superalgebra $D$ such that

$$
\mathcal{U}=\mathrm{var}_{\mathbb{Z}_2}^\ast(B)=\mathrm{var}_{\mathbb{Z}_2}^\ast(B_1\oplus
\cdots \oplus B_t \oplus D). \eqno (1)
$$

\noindent Moreover
$$
\mathrm{exp}_{\mathbb{Z}_2}^{\ast}(B_1)=\cdots
=\mathrm{exp}_{\mathbb{Z}_2}^{\ast}(B_t)=
\mathrm{exp}_{\mathbb{Z}_2}^{\ast}(\mathcal{U})=M+L
$$
 and
$$
\mathrm{exp}_{\mathbb{Z}_2}^{\ast}(D)<
\mathrm{exp}_{\mathbb{Z}_2}^{\ast}(\mathcal{U})= M+L.
$$
\smallskip

Let's now analyze the structure of a finite dimensional reduced $\ast$-superalgebra
 $R$  such that
$\mathrm{exp}_{\mathbb{Z}_2}^{\ast}(R)= M+L=
\mathrm{exp}_{\mathbb{Z}_2}^{\ast}(\mathcal{U})$ and
$\Gamma^\ast_{M^{\pm}+1, L^{\pm}+1}\subseteq
Id_{\mathbb{Z}_2}^{\ast}(R) $. We have that

$$
R=R_1\oplus \cdots \oplus R_m + J ,
 \eqno (2)$$

\noindent where $R_i$ are simple $\ast$-graded subalgebras of $R$,
$J=J(R)$ is the Jacobson radical of $R$ and $R_1J\cdots J R_m\neq
0$. By \cite[Theorem 4.3]{DSS} there exists a $\ast$-superalgebra
$\overline{R}$ isomorphic to the $\ast$-superalgebra
$UT_{\mathbb{Z}_2,\tilde{g}}^\ast(R_1, \ldots, R_m)$, for some
$\tilde{g}=(g_1,\ldots,g_m)\in \mathbb{Z}_2^m$, such that $Id(R)\subseteq Id(\overline{R})$ and

$$
\mathrm{exp}_{\mathbb{Z}_2}^{\ast}(R)=\mathrm{exp}_{\mathbb{Z}_2}^{\ast}(\overline{R})=\mathrm{exp}_{\mathbb{Z}_2}^{\ast}(UT_{\mathbb{Z}_2,\tilde{g}}^\ast(R_1,
\ldots, R_m)).
$$

\noindent It follows that
$$
M+L=\mathrm{exp}_{\mathbb{Z}_2}^{\ast}(R)=\mathrm{exp}_{\mathbb{Z}_2}^{\ast}(\overline{R})=
$$
$$
\mathrm{exp}_{\mathbb{Z}_2}^{\ast}(UT_{\mathbb{Z}_2,\tilde{g}}^\ast(R_1,
\ldots, R_m)) = \mathrm{dim}_FR_1+ \cdots + \mathrm{dim}_F R_m=
d_0^++d_0^-+d_1^+ + d_1^-
$$

\noindent where $d_i^{\pm}=\mathrm{dim}_F(R_1\oplus \cdots \oplus
R_m)_{(i)}^{\pm}$, for $i=0,1$.

Let  $0\leq \bar{m}\leq m$ denote the number of the
$\ast$-superalgebras $R_i$  with trivial grading appearing in (2).
We want to prove
that $\bar{m}=0$.

 Let's suppose  $\bar{m}> 0.$  By Lemma
\ref{capelliidentita}, $\overline{R}$ does not satisfy the $\ast$-graded
Capelli polynomials
$$ Cap_{d_0^++(m-\bar{m})+(\tilde{m}-1)+r_0}^{(\mathbb{Z}_2,
\ast)} [Y^+,X], \,\, \, \, Cap_{d_0^- +(m-\bar{m})+(\tilde{m}-1)+r_0}^{(\mathbb{Z}_2,
\ast)} [Y^-,X],
$$
$$Cap_{d_1^++(m-\bar{m})+(\tilde{m}-1)+r_1}^{(\mathbb{Z}_2,
\ast)} [Z^+,X], \,\,\,\, Cap_{d_1^-+(m-\bar{m})+(\tilde{m}-1)+r_1}^{(\mathbb{Z}_2,
\ast)} [Z^-,X],
$$
where $r_0$, $r_1$ are two non negative integers
dependent on the grading $\tilde{g}$ with $r_0+r_1=\bar{m}-\tilde{m}$.
 However  $\overline{R}$ satisfies $Cap_{M^++1}^{(\mathbb{Z}_2, \ast)}
[Y^+,X]$, $Cap_{M^-+1}^{(\mathbb{Z}_2, \ast)} [Y^-,X]$,
$Cap_{L^++1}^{(\mathbb{Z}_2, \ast)} [Z^+,X]$   and
$Cap_{L^-+1}^{(\mathbb{Z}_2, \ast)} [Z^-,X]$, then
$$
d_0^++(m-\bar{m})+(\tilde{m}-1)+r_0+d_0^-+(m-\bar{m})+(\tilde{m}-1)+r_0+
$$
$$
d_1^+ +(m-\bar{m})+(\tilde{m}-1)+r_1+ d_1^-+(m-\bar{m})+(\tilde{m}-1)+r_1\leq M+L.
$$

\smallskip
Since $d_0^++d_0^-+d_1^+ +d_1^-=M+L$ we obtain that $4(m-\bar{m})+4(\tilde{m}-1)+2(r_0+r_1)=0$  and so $2(m -1)+\tilde{m}-\bar{m}=0$
and this implies that $m\geq 2$. If $m=2$ then we easily obtain a contradiction. Thus $m=\bar{m}= \tilde{m}=1.$

Hence
$R=R_1\oplus J$ where $R_1 \simeq (M_{h_1}(F),t)$ or $R_1 \simeq
(M_{2h_1}(F),s)$ or $R_1 \simeq (M_{h_1}(F)\oplus M_{h_1}(F)^{op},
exc)$ with $h_1>0$.

 Now, let's analyze all possible cases as $M$ and
$L$ vary.
\smallskip

$\mathbf{1.}$ Let $M^{\pm}=\frac{h(h\pm 1)}{2}+\frac{l(l\pm 1)}{2}$ and $L^{\pm}=hl$, with $h\geq l > 0$.

 \noindent If $R \simeq (M_{h_1}(F),t)+J$ then
$\mathrm{exp}_{\mathbb{Z}_2}^{\ast}(R)=h_1^2$. Since
$\mathrm{exp}_{\mathbb{Z}_2}^{\ast}(R)=M+L=(h+l)^2$ we obtain that
$h_1=h+l$. By hypotesis, $R$ satisfies $Cap^{(\mathbb{Z}_2, \ast)}_{M^++1}[Y^+;X]$
but, since  $Id_{\mathbb{Z}_2}^\ast(R)\subseteq
Id_{\mathbb{Z}_2}^\ast(UT_{\mathbb{Z}_2,\tilde{g}}^\ast(R_1,\ldots,
R_q))$,  $R$ does not satisfy $Cap^{(\mathbb{Z}_2, \ast)} _{d_0^+}[Y^+;X].$ Hence,
for $h\geq l
>0$, we
have
      $$
      M^++1=\frac{h(h+1)}{2}+\frac{l(l+1)}{2}+1=\frac{h^2+l^2+(h+l)+2}{2}\leq
      $$
      $$
       \frac{h^2+l^2+(h+l)+2hl}{2}=\frac{(h+l)(h+l+1)}{2}=\frac{h_1(h_1+1)}{2}=d_0^+
      $$
      and this is impossible.

    \noindent   If $R \simeq (M_{2h_1}(F),s)+J$ then $\mathrm{exp}_{\mathbb{Z}_2}^{\ast}(R)=4h_1^2$. Since $\mathrm{exp}_{\mathbb{Z}_2}^{\ast}(R)=M+L=(h+l)^2$ we have that
    $2h_1=h+l$.
    Moreover $R$ satisfies
      $Cap^{(\mathbb{Z}_2, \ast)} _{M^-+1}[Y^-;X]$ but does not satisfy $Cap^{(\mathbb{Z}_2, \ast)} _{d_0^-}[Y^-;X]$ and
  so we get a contradiction since
      $$
      M^-+1=\frac{h(h-1)}{2}+\frac{l(l-1)}{2}+1=\frac{h^2+l^2-(h+l)+2}{2}<
      $$
      $$
       \frac{h^2+l^2+(h+l)+2hl}{2}=\frac{(h+l)^2+(h+l)}{2}=\frac{4h_1^2+2h_1}{2}=2h_1^2+h_1=d_0^-.
      $$

    \noindent   Finally, let $R \simeq (M_{h_1}(F)\oplus M_{h_1}(F)^{op}, exc)+J$, with $h_1>0$. Then $(h+l)^2=M+L=\mathrm{exp}_{\mathbb{Z}_2}^{\ast}(R)=2h_1^2,$  a
      contradiction.

      \noindent
  $\mathbf{2.}$ Let $M^{\pm}=h^2$ and $L^{\pm}=h(h\mp 1)$, with $h>0$.

  \noindent  If $R \simeq (M_{h_1}(F),t)+J$ then, as in the previous case, we obtain that $2h=h_1$. By hypothesis $R$
  satisfies $Cap^{(\mathbb{Z}_2, \ast)} _{M^++1}[Y^+;X]$ but it does not satisfy $Cap^{(\mathbb{Z}_2, \ast)} _{d_0^+}[Y^+;X]$, thus we have
  $$
  M^++1=h^2+1=(\frac{h_1}{2})^2+1=\frac{h_1^2}{4}+1\leq \frac{h_1^2}{2}+\frac{h_1}{2}=d_0^+
  $$
  a contradiction.

  \noindent If $R \simeq (M_{2h_1}(F),s)+J$ then  $h=h_1$. Since $R$ satisfies $Cap^{(\mathbb{Z}_2, \ast)} _{M^-+1}[Y^-;X]$ but does not satisfy
  $Cap^{(\mathbb{Z}_2, \ast)} _{d_0^-}[Y^-;X]$  we get the contradiction  $M^-+1=h^2+1=h_1^2+1<2h_1^2+h_1=d_0^-$.

  \noindent Finally, if $R \simeq (M_{h_1}(F)\oplus M_{h_1}(F)^{op}, exc)+J$
  with
  $h_1>0$, then we have $4h^2=2h_1^2,$ a contradiction.

  \noindent
  $\mathbf{3.}$ Let $M^{\pm}=h^2+l^2$ and $L^{\pm}=2hl$, with $h\geq l>0$.

  \noindent  If $R \simeq (M_{h_1}(F),t)+J$
   then we get the contradiction
    $2(h+l)^2=M+L=\mathrm{exp}_{\mathbb{Z}_2}^{\ast}(R)=h_1^2.$ The same occurs if  $R \simeq (M_{2h_1}(F),s)+J.$

 \noindent  Now, let $R \simeq (M_{h_1}(F)\oplus
  M_{h_1}(F)^{op},
  exc)+J$,
  with $h_1>0$. Then $2(h+l)^2=M+L=\mathrm{exp}_{\mathbb{Z}_2}^{\ast}(R)=2h_1^2$ and so $h_1=h+l$. Since $d_0^+=h_1^2$ we get that
  $M^++1=h^2+l^2+1<h^2+l^2+2hl=(h+l)^2=h_1^2=d_0^+$
  and this is
  impossible.

  \noindent
 $\mathbf{4., 5.}$ We consider the case $M^+=L^+=\frac{n(n+1)}{2}$ and $M^-=L^-=\frac{n(n-1)}{2}$. The proof of the other cases is very similar.

  \noindent If $R \simeq
 (M_{h_1}(F),t)+J$
 then $2n^2=M+L=\mathrm{exp}_{\mathbb{Z}_2}^{\ast}(R)=h_1^2,$ and
if $R \simeq (M_{2h_1}(F),s)+J$ then
 $2n^2=M+L=\mathrm{exp}_{\mathbb{Z}_2}^{\ast}(R)=4h_1^2$, a contradiction.

 \noindent Let
 $R \simeq (M_{h_1}(F)\oplus M_{h_1}(F)^{op}, exc)+J$, with $h_1>0$. Then $2n^2=M+L=\mathrm{exp}_{\mathbb{Z}_2}^{\ast}(R)=2h_1^2$ so $h_1=n$. Since $R$ satisfies
 $Cap^{(\mathbb{Z}_2, \ast)} _{M^-+1}[Y^-;X]$
 but it does not satisfy $Cap^{(\mathbb{Z}_2, \ast)} _{d_0^-}[Y^-;X]$ we have again a contradiction indeed $M^-+1=\frac{n(n-1)}{2}+1\leq n(n-1)+1\leq n^2=h_1^2=d_0^-$.

 \noindent
 $\mathbf{6.}$ Let $M^{\pm}=L^{\pm}=n^2$,  with $n>0$.

   \noindent If $R \simeq (M_{h_1}(F)\oplus M_{h_1}(F)^{op}, exc)+J$ then $4n^2=M+L=\mathrm{exp}_{\mathbb{Z}_2}^{\ast}(R)=2h_1^2$ a
 contradiction.

 \noindent  If
 $R \simeq (M_{h_1}(F),t)+J$ then $4n^2=M+L=\mathrm{exp}_{\mathbb{Z}_2}^{\ast}(R)=h_1^2$ and so $h_1=2n$.
 $R$ satisfies $Cap^{(\mathbb{Z}_2, \ast)} _{M^++1}[Y^+;X]$ but does not
 satisfy
 $Cap^{\ast}_{d_0^+}[Y^+;X]$  then we obtain a contradiction in fact
 $M^++1=n^2+1=\frac{h_1^2}{4}+1 \leq \frac{h_1^2}{2}+ \frac{h_1}{2}=\frac{h_1(h_1+1)}{2}=d_0^+$.

 \noindent  Finally, let $R
 \simeq (M_{2h_1}(F),s)+J$. Hence $4n^2=M+L=\mathrm{exp}_{\mathbb{Z}_2}^{\ast}(R)=4h_1^2$ and so $n=h_1$. Also in this case we get the contradiction
 $M^-+1=n^2+1<2n^2+1<2n^2+n=2h_1^2+h_1=d_0^-$.

\medskip
So we obtained that $\bar{m}=0.$

\smallskip

Let   $ R=R_1\oplus \cdots \oplus R_m + J, $ where $R_i$ are simple
$\ast$-superalgebras with non trivial grading.  Let's prove that $m
= 1.$
  By Lemma \ref{capelliidentita}, $\overline{R}$ does not satisfy
 the $\ast$-graded Capelli polynomials
$Cap_{d_0^++m-1}^{(\mathbb{Z}_2, \ast)} [Y^+,X]$, $Cap_{d_0^-
+m-1}^{(\mathbb{Z}_2, \ast)} [Y^-,X]$,
$Cap_{d_1^++m-1}^{(\mathbb{Z}_2, \ast)} [Z^+,X]$   and
$Cap_{d_1^-+m-1}^{(\mathbb{Z}_2, \ast)} [Z^-,X]$
 but  satisfies $Cap_{M^++1}^{(\mathbb{Z}_2, \ast)} [Y^+,X]$, $Cap_{M^-+1}^{(\mathbb{Z}_2, \ast)} [Y^-,X]$,
$Cap_{L^++1}^{(\mathbb{Z}_2, \ast)} [Z^+,X]$   and $Cap_{L^-+1}^{(\mathbb{Z}_2, \ast)} [Z^-,X]$ thus $d_0^++m-1\leq M^+$, $d_0^-+m-1\leq M^-$, $d_1^++m-1\leq L^+$ and
$d_1^-+m-1\leq L^-$. Hence we have that
$$
d_0^++(m-1)+d_0^-+(m-1)+d_1^++(m-1)+d_1^-+(m-1)\leq M^++M^-+L^++L^-=M+L.
$$
Since $d_0^++d_0^-+d_1^-+d_1^-=M+L$ we obtain that $4(m-1)=0$ and so $m=1$.

\smallskip It follows that $R=R_1\oplus J$ where $R_1$ is a simple
$\ast$-superalgebra with non trivial grading. Now let's analyze the
cases corresponding to the different values of M and L.

\smallskip

\noindent
$\mathbf{1.}$ Let $M^{\pm}=\frac{h(h\pm 1)}{2}+\frac{l(l\pm 1)}{2}$ and $L^{\pm}=hl$, with $h\geq l > 0$.

 \noindent If $R \simeq (M_{h_1, h_1}(F),s)+J$ then
$(h+l)^2=M+L=\mathrm{exp}_{\mathbb{Z}_2}^{\ast}(R)=4h_1^2$ so we
have $2h_1=h+l$. By hypothesis $R$ satisfies $Cap^{(\mathbb{Z}_2,
\ast)}_{L^-+1}[Z^-;X]$ but does not satisfy $Cap^{(\mathbb{Z}_2,
\ast)}_{d_1^-}[Z^+;X]$, where $d_1^-=h_1(h_1+1)$. Since $h+l=2h_1$
and $h\geq l > 0$ we have that $h_1^2\geq hl$ and so

     $$
     L^-+1=hl+1\leq h_1^2+1\leq h_1(h_1+1)=d_1^-
     $$
     a contradiction.

  \noindent    If $R \simeq (M_{h_1, l_1}(F)\oplus M_{h_1,l_1}(F)^{op}, exc)+J$, with $h_1\geq l_1 >0$,
     then $(h+l)^2=M+L=\mathrm{exp}_{\mathbb{Z}_2}^{\ast}(R)=2(h_1+l_1)^2$ and so
     we
     have  again a
     contradiction.

   \noindent   If $R \simeq (M_n(F+cF),\ast)+J$, where $(a+cb)^\ast= a^\diamond \pm cb^\diamond$ and $\diamond =t,s$, then we obtain the contradiction $(h+l)^2=2n^2$.

     \noindent If  $R \simeq (M_n(F+cF)\oplus M_n(F+cF)^{op}, exc)+J$ with $n>0$, then $(h+l)^2=M+L=\mathrm{exp}_{\mathbb{Z}_2}^{\ast}(R)=4n^2$ and so $2n=h+l$. As before we
     can easily
     obtain a
     contradiction. It follows that $R \simeq (M_{h, l}(F),t)+J$.
\smallskip

 \noindent
  $\mathbf{2.}$ Let now $M^{\pm}=h^2$ and $L^{\pm}=h(h\mp 1)$, with $h>0$.

   \noindent If $R \simeq (M_{h_1, l_1}(F),t)+J$, then, since $M+L=\mathrm{exp}_{\mathbb{Z}_2}^{\ast}(R)$, we have
  $4h^2=(h_1+l_1)^2$
  and so $h_1+l_1=2h^2$. By hypothesis $R$ satisfies $Cap^{(\mathbb{Z}_2, \ast)}_{M^++1}[Y^+;X]$ but does not satisfy $Cap^{(\mathbb{Z}_2, \ast)}_{d_0^+}[Y^+;X]$  where
  $d_0^+=\frac{h_1(h_1+1)}{2}+\frac{l_1(l_1+1)}{2}$. Since $h_1+l_1=2h$ and $h_1\geq l_1 >0$ we have $h^2\geq h_1l_1$ and so it follows that
     $$
     M^++1=h^2+1<h(2h+1)-h_1l_1=\frac{h_1+l_1}{2}(h_1+l_1+1)-h_1l_1=
     $$
     $$
     \frac{h_1(h_1+1)}{2}+\frac{l_1(l_1+1)}{2}=d_0^+
     $$
     a contradiction.

     \noindent  If $R \simeq (M_{h_1, l_1}(F)\oplus M_{h_1,l_1}(F)^{op}, exc)+J$, with $h_1\geq l_1 >0$, or $R \simeq (M_n(F+cF),\ast)+J$
       where $(a+cb)^\ast= a^\diamond \pm cb^\diamond$ and $\diamond =t,s$ then
    easily we get a contradiction.

 \noindent If  $R \simeq (M_n(F+cF)\oplus M_n(F+cF)^{op}, exc)+J$ with $n>0$, then $4h^2=M+L=\mathrm{exp}_{\mathbb{Z}_2}^{\ast}(R)=4n^2$ and so $n=h$.
   $R$  satisfies
     $Cap^{(\mathbb{Z}_2, \ast)}_{L^-+1}[Z^-;X]$ but $R$ does not satisfy $Cap^{(\mathbb{Z}_2, \ast)}_{d_1^-}[Z^+;X]$, where $d_1^-=n^2$ and  we obtain the following
     contradiction
     $
     L^-+1=h(h-1)=h^2-h-1\leq h^2=n^2=d_1^-.
     $
 So, in this case,  $R \simeq (M_{h, h}(F),s)+J$.

\smallskip \noindent
  $\mathbf{3.}$  Let $M^{\pm}=h^2+l^2$ and $L^{\pm}=2hl$, with $h\geq l>0$.

   \noindent  If $R \simeq (M_{h_1, l_1}(F),t)+J$, $R \simeq  (M_{h_1, h_1}(F),s)+J$ or  $R \simeq
  (M_n(F+cF)\oplus
  M_n(F+cF)^{op}, exc)+J$ easily  we get  a contradiction.

 \noindent   If $R \simeq (M_n(F+cF),\ast)+J$ where $(a+cb)^\ast= a^\diamond \pm cb^\diamond$ and $\diamond =t,s$ then we
  have that $2(h+l)^2=2n^2$ and so $h+l=n$. Let consider the case when $R \simeq (M_n(F+cF),\ast)+J$ with $(a+cb)^\ast= a^t - cb^t$, the other cases are very similar. Since
  $R$
  satisfies $Cap^{(\mathbb{Z}_2, \ast)}_{L^-+1}[Z^-;X]$ but $R$ does not satisfy $Cap^{(\mathbb{Z}_2, \ast)}_{d_1^-}[Z^+;X]$ we obtain
      $$
      L^-+1= 2hl+1<\frac{(h+l+1)(h+l)}{2}=\frac{(n+1)n}{2}=d_1^-
      $$
      a contradiction. It follows that $R \simeq (M_{h, l}(F)\oplus M_{h_1,l_1}(F)^{op}, exc)+J$.

 \smallskip\noindent
  $\mathbf{4., 5.}$  Let consider the case $M^+=L^+=\frac{n(n+1)}{2}$ and $M^-=L^-=\frac{n(n-1)}{2}$. The proof of the other cases is very similar. As before let $R \simeq
  (M_n(F+cF)\oplus M_n(F+cF)^{op}, exc)+J$, then $2n^2=2(h_1+l_1)^2$ and so $n=h_1+l_1$ with $h_1\geq l_1>0$.
   $R$ satisfies $Cap^{(\mathbb{Z}_2, \ast)}_{M^-+1}[Y^-;X]$ but does not
  satisfy
  $Cap^{(\mathbb{Z}_2, \ast)}_{d_0^-}[Y^-;X]$ then
      $$
      \frac{n(n-1)}{2}+1=\frac{n^2-n+2}{2}\leq \frac{n^2-1}{2}=\frac{(h_1+l_1)^2-1}{2}<h_1^2+l_1^2=d_0^-
      $$
      a contradiction.
 In all other cases we obtain a contradiction except when $R \simeq (M_n(F+cF),\ast)+J$ and $(a+cb)^\ast= a^t + cb^t$.

  \smallskip\noindent
  $\mathbf{6.}$ Let $M^{\pm}=L^{\pm}=n^2$,  with $n>0$.

   \noindent If $R \simeq (M_{h_1, l_1}(F)\oplus M_{h_1, l_1}(F)^{op}, exc)+J$ or $R \simeq (M_n(F+cF),\ast)+J$ with $(a+cb)^\ast=
  a^\diamond \pm cb^\diamond$ and $\diamond =t,s$, then easily  we get  a contradiction.

   \noindent If $R \simeq (M_{h_1, l_1}(F),t)+J$, then $h_1+l_1=2n$ and
      with  analogous reasoning to that of case 2 we obtain a contradiction.

     \noindent  So let assume that $R \simeq (M_{h_1, h_1}(F),s)+J$, then $4n^2=4h_1^2$ and so $h_1=n$. Because $R$ satisfies  $Cap^{(\mathbb{Z}_2, \ast)}_{L^-+1}[Z^-;X]$ but
     it does not satisfy
      $Cap^{(\mathbb{Z}_2, \ast)}_{d_1^-}[Z^-;X]$ we obtain
      $L^-+1=n^2+1\leq n(n+1) = h_1(h_1+1)=d_1^-$ and this is impossible.
It follows that $R \simeq (M_n(F+cF)\oplus M_n(F+cF)^{op}, exc)+J$.

\smallskip
Thus we have proved that $R \simeq A+J$ where $A$ is a simple
$\ast$-superalgebra with non trivial grading. Then, from Lemmas
\ref{radicaltrasposta}, \ref{j10trasposta}, \ref{Ntrasposta} we
obtain that

 $$
R\cong (A+J_{11})\oplus J_{00}\cong (A\otimes
N^\sharp)\oplus J_{00}
$$

\noindent where $N^\sharp $ is the algebra obtained from $N$ by
adjoining a unit element. Since $N^\sharp $ is commutative, it
follows that $A + J_{11}$ and $A$ satisfy the same $\ast$-graded
identities. Thus
$\mathrm{var}_{\mathbb{Z}_2}^\ast(R)=\mathrm{var}_{\mathbb{Z}_2}^\ast(A\oplus
J_{00})$ with $J_{00}$ finite dimensional nilpotent
$\ast$-superalgebra. Hence, from the decomposition $(1)$, we get

$$
\mathcal{U}=\mathrm{var}_{\mathbb{Z}_2}^\ast(\Gamma^\ast_{M^{\pm}+1,L^{\pm}+1})=\mathrm{var}_{\mathbb{Z}_2}^{\ast}(A\oplus
D),
$$

\noindent where $D$ is a finite dimensional $\ast$-superalgebra with
$\mathrm{exp}_{\mathbb{Z}_2}^{\ast}(D)<M+L$ and the theorem is proved.

\bigskip

From Corollary \ref {supercodimensioni} we easily obtain the following

\bigskip

\begin{cor}
\begin{enumerate}
\item[1)] If $M^{\pm}=\frac{h(h\pm 1)}{2}+\frac{l(l\pm 1)}{2}$ and $L^{\pm}=hl$, with $h\geq l > 0$, then
$$
c_n^{(\mathbb{Z}_2,
\ast)}(\Gamma^\ast_{M^{\pm}+1,L^{\pm}+1})
\simeq
c^{(\mathbb{Z}_2,\ast)}_n((M_{h,l}(F),t));
$$
  \item[2)] If $M^{\pm}=h^2$ and $L^{\pm}=h(h\mp 1)$, with $h>0$, then
$$
c_n^{(\mathbb{Z}_2,
\ast)}(\Gamma^\ast_{M^{\pm}+1,L^{\pm}+1})\simeq
c^{(\mathbb{Z}_2,\ast)}_n((M_{h,h}(F),s));
$$
\item[3)] If $M^{\pm}=h^2+l^2$ and $L^{\pm}=2hl$, with $h\geq l > 0$, then
$$
c_n^{(\mathbb{Z}_2,
\ast)}(\Gamma^\ast_{M^{\pm}+1,L^{\pm}+1})\simeq
c^{(\mathbb{Z}_2,\ast)}_n((M_{h,l}(F)\oplus M_{h,l}(F)^{op},exc));
$$
 \item[4)] If $M^{+}=L^{\pm}=\frac{n(n+1)}{2}$, $M^{-}=L^{\mp}=\frac{n(n-1)}{2}$,  with $n>0$, then
$$
c_n^{(\mathbb{Z}_2,\ast)}(\Gamma^\ast_{M^{\pm}+1,L^{\pm}+1})\simeq
c^{(\mathbb{Z}_2,\ast)}_n((M_n(F+cF),\ast))
$$
where $(a+cb)^\ast= a^t\pm cb^t;$
\item[5)]  If $M^{+}=L^{\pm}=\frac{n(n-1)}{2}$, $M^{-}=L^{\mp}=\frac{n(n+1)}{2}$,  with $n>0$, then
$$
c_n^{(\mathbb{Z}_2,
\ast)}(\Gamma^\ast_{M^{\pm}+1,L^{\pm}+1})\simeq
c^{(\mathbb{Z}_2,\ast)}_n((M_n(F+cF),\ast)),
$$
 where $(a+cb)^\ast= a^s
\pm cb^s$;
\item [6)] If $M^{\pm}=L^{\pm}=n^2$,  with $n>0$,  then
$$
c_n^{(\mathbb{Z}_2,
\ast)}(\Gamma^\ast_{M^{\pm}+1,L^{\pm}+1})\simeq
c^{(\mathbb{Z}_2,\ast)}_n(M_n(F+cF)\oplus M_n(F+cF)^{op},exc)).
$$
\end{enumerate}

\end{cor}

\bigskip

\bibliographystyle{amsplain}

\end{document}